\documentclass[11pt,openbib]{article} 
\usepackage{amsfonts,amsthm}
\usepackage{amssymb,amsmath}
\newcommand{\R}{\mathbb{R}} 
\newcommand{\C}{\mathbb{C}}
 
\newcommand{\Z}{\mathbb{Z}}
 
\newcommand{\gl}{\mathop{\mathrm{gl}}\nolimits}

\newcommand{\ssl}{\mathop{\mathrm{sl}}\nolimits} 
\newcommand{\ad}{\mathop{\mathrm{ad}}\nolimits}

\newcommand{\comp}{\raisebox{0pt}{$\scriptstyle\circ \, $}}

\newcommand{\setrule}{\, \mathop{\rule[-4pt]{.5pt}{13pt}\, }\nolimits}
\newcommand{\smallspace}{\smallskip\par\noindent}

\newcommand{\rowspace}{\rule{0pt}{16pt}}

\newcommand{\onehalf}{\mbox{$\frac{\scriptstyle 1}{\scriptstyle 2}\,$}} 
\newcommand{\spann}{\mathop{\rm span}\nolimits} 
\newcommand{\ttfrac}[2]{\mbox{$\frac{{\scriptstyle #1}}{{\scriptstyle #2}}$}}
\allowdisplaybreaks

\begin{document}
\begin{center}
{\Large \bf System of roots}
\footnotetext{version: \today} \\ 
\vspace{.05in} 
Richard Cushman 
\end{center}

In this paper we define the notion of a system of roots, which generalizes the concept of 
root system of a semisimple complex Lie algebra \cite[p.42]{humphreys}. We show 
every algebra in a certain class of nonsemisimple Lie algebras has a system of roots. This leads to an analogue of a Weyl group, which we study in another paper \cite{cushman17b}. 

\section{System of roots}
\label{sec1} 
 
In this section we give the axioms of a system of roots, which unlike the axioms for a root 
system does \emph{not} use an inner product. We prove some consequences of these 
axioms, one of which is that every root system is a system of roots.   

\subsection{Axioms for a system of roots}
\label{sec1sub1}

Let $V$ be a finite dimensional real vector space with $\Phi $ a finite subset of nonzero vectors, 
which satisfy the following axioms. \medskip 

\noindent 1. $V = {\spann }_{\R }\Phi $, using addition $+$ of vectors in $V$.
\smallspace 
\noindent 2. $\Phi = - \Phi $, where $-$ is the additive inverse of $+$. 
\smallspace
\noindent 3. \parbox[t]{4.75in}{For every $\beta $, $\alpha \in \Phi $ there is an \emph{extremal root chain 
${\mathcal{S}}^{\beta }_{\alpha }$ through $\beta $ in the direction $\alpha $} given by 
$\{ \beta +j \alpha \in \Phi \cup \{ 0 \} \setrule \mbox{for every $j \in \Z$, $-q \le j \le p $} \} $. Here 
$q,p \in {\Z }_{\ge 0}$ and are as large as possible. The pair $(q,p)$ is the \linebreak 
\emph{integer pair} associated to ${\mathcal{S}}^{\beta }_{\alpha }$. The integer $\langle \beta , \alpha \rangle = q-p$ 
is called the \emph{Killing integer} of ${\mathcal{S}}^{\beta }_{\alpha }$.} 
\smallspace
\noindent  4. Fix $\alpha \in \Phi $ and suppose that ${\beta }_1$, ${\beta }_2$, and 
${\beta }_1+ {\beta }_2 \in \Phi $. Then 
\begin{equation}
\langle {\beta }_1+{\beta }_2, \alpha \rangle = \langle {\beta }_1, \alpha \rangle + \langle {\beta }_2, \alpha \rangle .
\label{eq-sec1subsec1one}
\end{equation}
\noindent 5. For every $\alpha \in \Phi $ we have $\langle \alpha , \alpha \rangle = 2$. \bigskip 

\noindent We call $\Phi \cup \{ 0 \} $ a \emph{system of roots} in $V$. \medskip 

The main point in the axioms for a system of roots is that there is \emph{no} Euclidean inner 
product on $V$. In fact, this distinguishes a system of roots from a root system of a semisimple complex 
Lie algebra.

\subsection{Some consequences}
\label{sec1subsec2}

In this section we draw some consequences from the axioms of a system of roots. \medskip 

\noindent \textbf{Lemma 1.2.1} Axiom $3$ is follows from axiom $1$. \medskip 

\noindent \textbf{Proof.} Axiom 3 holds because for every $\beta $, $\alpha \in \Phi $ the affine line
$\R \rightarrow V: t \mapsto \beta +t \alpha $ through $\beta $ in the direction of $\alpha $ intersects 
$\Phi \cup \{ 0 \} $ in at most a finite number of distinct values of the parameter $t \in \Z $. Since 
$\beta \in \Phi $ it follows that $0$ is such a parameter value. Thus there is a maximal $q$ and $p 
\in {\Z }_{\ge 0}$ such that for every $j \in \Z $ with $-	q \le j \le p$ we have $\beta + j \alpha \in 
\Phi \cup \{ 0 \} $. So an extremal root chain through $\beta $ in the direction $\alpha $ exists. By 
definition $\langle \beta , \alpha \rangle = q-p$. Thus axiom $3$ is a consequence of axiom $1$. \medskip 

\noindent \textbf{Claim 1.2.2} Let $\alpha $, $\beta \in \Phi $. If $\langle \beta , \alpha \rangle < 0$, then 
$\beta + \alpha \in \Phi \cup \{ 0 \}$; while if $\langle \beta , \alpha \rangle > 0 $, then $\beta - \alpha 
\in \Phi \cup \{ 0 \} $. \medskip

\noindent \textbf{Proof.} From axiom $3$ we deduce that the extremal root chain ${\mathcal{S}}^{\beta }_{\alpha }$ 
through $\beta $ in the direction $\alpha $ with integer pair $(q,p)$ exists. If $q-p = \langle \beta , \alpha \rangle 
< 0$, then $p > q \ge 0$. So ${\mathcal{S}}^{\beta }_{\alpha }$ contains $\beta + \alpha $. Therefore 
$\beta + \alpha \in \Phi \cup \{ 0 \} $. If $q -p = \langle \beta , \alpha \rangle > 0$, then $q > p \ge 0$. So 
the extremal root chain ${\mathcal{S}}^{\beta }_{\alpha }$ through $\beta $ in the direction $\alpha $ contains 
$\beta - \alpha $. Therefore $\beta - \alpha \in \Phi \cup \{ 0 \} $. \hfill $\square $ \medskip 

\noindent \textbf{Lemma 1.2.3} For every $\beta $, $\alpha \in \Phi $ we have $\langle \beta , - \alpha \rangle 
= - \langle \beta , \alpha \rangle $ and $\langle - \beta , \alpha \rangle = - \langle \beta , \alpha \rangle $. \medskip 

\noindent \textbf{Proof.} From axiom $3$ the extremal root chain ${\mathcal{S}}^{\beta }_{\alpha }: 
\beta -q \alpha , \ldots , \beta + p \alpha $ through $\beta $ in the direction $\alpha $ with integer pair 
$(q,p)$ exists. From axiom $2$ it follows that $-\beta $ and $-\alpha \in \Phi $. The following argument 
shows that the extremal root chain ${\mathcal{S}}^{\beta }_{-\alpha }$ has integer pair $(p,q)$. From 
the extremal root chain ${\mathcal{S}}^{\beta }_{\alpha }$ we see that 
\begin{displaymath}
\beta -p(-\alpha ), \, \beta -(p-1)(-\alpha ) , \ldots , \beta +(q-1)(-\alpha ) , \,  \beta + q(-\alpha )
\end{displaymath}
is an extremal root chain through $\beta $ in the direction $-\alpha $ with integer pair $(p,q)$. In other words, 
$\langle \beta , -\alpha \rangle = p -q = - \langle \beta , \alpha \rangle $. 
\par Next we show that the extremal root chain ${\mathcal{S}}^{-\beta }_{\alpha }$ through $-\beta $ in the 
direction $\alpha $ has integer pair $(p,q)$. Multiplying the elements of the extremal root chain 
${\mathcal{S}}^{\beta }_{\alpha }$ by $-1$ and using axiom $2$, we see that 
\begin{displaymath}
-\beta - p \alpha , \, -\beta -(p-1)\alpha , \ldots , -\beta +(q-1)\alpha , \, -\beta + q\alpha 
\end{displaymath}
is an extremal root chain ${\mathcal{S}}^{-\beta }_{\alpha }$ through $-\beta $ in the direction $\alpha $ 
with integer pair $(p,q)$. In other words, $\langle - \beta , \alpha \rangle = p-q = - \langle \beta , \alpha \rangle $. \hfill $\square $ \medskip 

We now consider root systems, see \cite[p.42]{humphreys}. Let $\big( U, (\, \, , \, \, ) \big) $ be a 
finite dimensional real vector space with a Euclidean inner product $(\, \, , \, \, ) $. Let $\Phi $ be a 
finite subset of nonzero vectors in $U$. Suppose that the following axioms hold. \medskip 

\noindent 1. $U = {\spann }_{\R }\Phi $. \smallspace
2. If $\alpha \in \Phi $ and $\lambda \in \R $, then $\lambda \alpha \in \Phi $ if and only if $|\lambda | =1$. 
\smallspace
3. \parbox[t]{4.25in}{If $\alpha \in \Phi $, then the reflection ${\sigma }_{\alpha }: U \rightarrow U: v \mapsto v - 
\ttfrac{2(v,\alpha )}{(\alpha , \alpha )} \alpha $ is an orthogonal real linear mapping of $\big( U, (\, \, , \, \, ) \big)$ into itself, which preserves $\Phi $.} 
\smallspace
4. If $\beta $, $\alpha \in \Phi $, then $\langle \beta , \alpha \rangle = \frac{2(\beta ,\alpha )}{(\alpha , \alpha )} 
\in \Z $. \medskip 

\noindent Then $\Phi $ is a \emph{root system}. \medskip 

\noindent \textbf{Claim 1.2.4} Every root system is a system of roots. \medskip 

\noindent \textbf{Proof.} Suppose that $\Phi $ is a root system. Then axiom $1$ of root system is the 
same as axiom $1$ of a system of roots. If $\alpha \in \Phi $, then from axiom $2$ of a root system 
it follows that $-\alpha \in \Phi $ and thus $\alpha = -(-\alpha )$. So $\Phi = - \Phi $, which is axiom $2$ 
of a system of roots. From lemma 1.2.1 it follows that for every $\beta $, $\alpha \in \Phi $ there is an 
extremal root chain ${\mathcal{S}}^{\beta }_{\alpha }$ through $\beta $ in the direction $\alpha $ with 
integer pair $(q,p)$. Thus the first statement in axiom $3$ for a system of roots holds. This does not 
complete its verification, because we still have to show that $\langle \beta , \alpha \rangle = q-p$, where 
$\langle \beta , \alpha \rangle = \frac{2(\beta , \alpha )}{(\alpha , \alpha )}$. Axiom $5$ of a system 
of roots follows because for every $\alpha \in \Phi $, which is nonzero by hypothesis, we have 
$\langle \alpha , \alpha \rangle = \frac{2(\alpha , \alpha )}{(\alpha , \alpha )} = 2$. From the 
definition $\langle \beta , \alpha \rangle = \frac{2(\beta , \alpha )}{(\alpha , \alpha )}$ it follows that 
for each fixed $\alpha \in \Phi $ the function 
$K_{\alpha }: \Phi \rightarrow \Z : \beta \mapsto \langle \beta , \alpha \rangle $ is linear and is $\Z $-valued 
by axiom $4$ of a root system. This proves axiom 4 of a system of roots. \medskip 

We now finish the proof of axiom $3$ of a system of roots. Using axiom $3$ of a root system we show 
that for every $\alpha \in \Phi $ the orthogonal reflection ${\sigma }_{\alpha }:U \rightarrow U: 
v \mapsto v - K_{\alpha }(v)\alpha $ maps the extremal root chain ${\mathcal{S}}^{\beta }_{\alpha }$ 
into itself. For every $j \in \Z $ with $-q \le j \le p$ let $\beta +j \alpha \in {\mathcal{S}}^{\beta }_{\alpha }$. Then 
\begin{align}
{\sigma }_{\alpha }(\beta +j \alpha ) & = \beta +j \alpha -K_{\alpha }(\beta +j \alpha ) \alpha \notag \\
& = \beta +j \alpha -K_{\alpha }(\beta )\alpha - jK_{\alpha }(\alpha ) \alpha \notag \\
& = \beta - (\langle \beta , \alpha \rangle +j) \alpha . \notag 
\end{align}  
Since $\beta +j \alpha \in \Phi \cup \{ 0 \} $ and ${\sigma }_{\alpha }$ maps $\Phi $ into itself and $0$ into $0$, 
it follows that $\beta - (\langle \beta , \alpha \rangle +j) \alpha \in \Phi \cup \{ 0 \} $. Because 
$\langle \beta , \alpha \rangle \in \Z $, for every $j \in \Z$ with $-q \le j \le p$ we have 
$-(\langle \beta , \alpha \rangle +j ) \in F = \big\{ \ell \in \Z, \, -q \le \ell \le p \setrule \beta + \ell \alpha \in \Phi \cup \{ 0 \} \big\} $. So the orthogonal reflection ${\sigma }_{\alpha }$ maps the extremal root chain 
${\mathcal{S}}^{\beta }_{\alpha }$ 
into itself. Moreover $-(\langle \beta , \alpha \rangle +p) $ is the smallest element of $F$. But 
${\mathcal{S}}^{\beta }_{\alpha }$ is an extremal root chain with integer pair $(q,p)$. So 
$-q = -\langle \beta , \alpha \rangle -p$, that is, $\langle \beta , \alpha \rangle = q-p$. This proves axiom $3$ of 
a system of roots and thereby the claim. \hfill $\square $ 

\section{Fundamental sandwich algebras}
\label{sec2}

The goal of this section is to show that every fundamental sandwich algebra 
has a system of roots. In the appendix we list the fundamental sandwich algebras. \medskip 

We start by defining a fundamental sandwich algebra. A complex Lie algebra 
$\widetilde{\mathfrak{g}}$ is a \emph{sandwich algebra} if it satisfies the two conditions: 
1) $\widetilde{\mathfrak{g}}$ is the direct sum of a simple complex Lie algebra $\mathfrak{g}$ 
and a nilpotent ideal $\widetilde{\mathfrak{n}}$, which is a sandwich, namely, $[\widetilde{\mathfrak{n}}, [\widetilde{\mathfrak{n}}, \widetilde{\mathfrak{n}}]] =0$ and $[\widetilde{\mathfrak{n}}, \widetilde{\mathfrak{n}}] 
\ne 0$.\footnote{Since $\widetilde{\mathfrak{g}}$ is the 
nilradical of $\widetilde{\mathfrak{g}}$, the Lie algebra $\widetilde{\mathfrak{g}}$ is not semisimple.} 
2) Let $\mathfrak{h}$ be a Cartan subalgebra of $\mathfrak{g}$ with associated root system 
$\mathcal{R}$. Then $\widetilde{\mathfrak{n}}$ is an ${\ad }_{\mathfrak{h}}$-invariant subspace of $\widetilde{\mathfrak{g}}$, because $\widetilde{\mathfrak{n}}$ is an ideal. We assume that 
${\ad }_{\mathfrak{h}}|\widetilde{\mathfrak{n}}$ is a maximal toral subalgebra of 
$\gl (\widetilde{\mathfrak{n}}, \C )$. The set of roots $\mathfrak{R}$ of 
${\ad }_{\mathfrak{h}}|\widetilde{\mathfrak{n}}$ is a set of $\tau \in {\mathfrak{h}}^{\ast }$ 
such that  for some nonzero $X_{\tau } \in \widetilde{\mathfrak{n}}$ we have 
$[H, X_{\tau }] = \tau (H) X_{\tau } $ for all $H \in \mathfrak{h}$. The sandwich algebra 
$\widetilde{\mathfrak{g}}$ is \emph{very special} if it satisfies the following three conditions: 1) there is  complex simple Lie algebra $\underline{\mathfrak{g}}$, having $\widetilde{\mathfrak{g}}$ as a subalgebra, which has rank $1$ greater than the rank of $\mathfrak{g}$ and $\widetilde{\mathfrak{n}}$ is an 
ideal of $\underline{\mathfrak{g}}$. 2) there is a Cartan subalgebra 
$\underline{\mathfrak{h}}$ of $\underline{\mathfrak{g}}$ with associated root system 
$\underline{\mathcal{R}}$, which is \emph{aligned} with the Cartan subalgebra 
$\mathfrak{h}$ of $\widetilde{\mathfrak{g}}$, that is, there is a nonzero 
vector\footnote{Since $\underline{\mathfrak{g}}$ is a simple Lie algebra $\underline{E} = 
{\spann}_{\R}\underline{\mathcal{R}}$ is a real vector subspace of ${\underline{\mathfrak{h}}}^{\ast }$, 
whose complexification is ${\underline{\mathfrak{h}}}^{\ast }$, such that the Killing form 
$\underline{k}$ on $\underline{\mathfrak{g}}$ induced on $\underline{E}$ is the Euclidean 
inner product, see \cite[p.39--40]{humphreys}. The map ${\underline{k}}^{\sharp}: 
\underline{E} \rightarrow {\underline{E}}^{\ast }: \underline{\alpha } \mapsto 
{\underline{k}}^{\sharp}(\underline{\alpha})$, where 
${\underline{k}}^{\sharp}(\underline{\alpha}){\underline{\alpha}}\, ' = 
\underline{k}(\underline{\alpha }, {\underline{\alpha }}\, ')$, is an isomorphism with inverse 
${\underline{k}}^{\flat}$.} $\underline{\widetilde{H}} \in 
{\underline{k}}^{\flat}(\underline{E}) \subseteq \underline{\mathfrak{h}}$ such that 
${\mathfrak{h}}^{\ast } = \{ \underline{\alpha } \in \underline{\mathcal{R}} \setrule \, 
\underline{\alpha }(\underline{\widetilde{H}}) = 0 \} $. \linebreak 
3) $\widehat{\alpha } \in \mathcal{R}$ if and 
only if for every $\underline{\alpha } \in \underline{\mathcal{R}}$ such that $\underline{\alpha}|\mathfrak{h} = 
\widehat{\alpha}$ we have $\underline{\alpha }(\underline{\widetilde{H}}) < 0$. Let 
$E = {\spann }_{\R }\mathcal{R}$. Since ${\underline{k}}^{\flat}(\underline{E}) = {\underline{k}}^{\flat}(E) \oplus 
{\spann }_{\R } \{ \underline{\widetilde{H}} \} $ and $\underline{E}$ is finite dimensional, there is 
a unique $\underline{\alpha} \in \underline{\mathcal{R}}$ satisfying condition 3). 
A very special sandwich algebra is \emph{fundamental} if the center of the nilradical 
$\widetilde{\mathfrak{n}}$ is spanned by the vector $X_{\zeta }$ where $\zeta = 
-(\underline{\widetilde{H}})^{\ast }|\mathfrak{h}$. Here $(\underline{\widetilde{H}})^{\ast }$ is the 
vector ${\underline{k}}^{\sharp}(\underline{\widetilde{H}}) \in \underline{E}$. 
Since $\underline{k}$ is the Euclidean inner product on $\underline{E}$, we have 
$\underline{k}\big( {\underline{k}}^{\sharp}(\underline{\widetilde{H}}), 
{\underline{k}}^{\sharp}(\underline{\widetilde{H}}) \big) >0$. 
So $\zeta (\underline{\widetilde{H}} ) < 0$. Thus $\zeta \in \mathfrak{R}$, which 
implies $\zeta (H) =0$ for all $H \in \mathfrak{h}$, since the Cartan subalgebras 
$\mathfrak{h}$ and $\underline{\mathfrak{h}}$ are aligned. Hence $\zeta $ is the zero root 
$\widehat{0}$ in $\mathfrak{R}$. \medskip 

We now give an outline of the argument in this section which shows that every fundamental sandwich algebra has a system of roots. In \S 2.1 we define an operation of addition on 
$\widehat{\mathcal{R}} = \mathfrak{R} \setminus \{ \zeta \}$. This allows us to show that axioms 1 and 2 for a system of roots holds for the set $\Phi = \widehat{\mathcal{R}}$. Next we construct a set of positive simple roots $\widehat{\Pi}$ for $\widehat{\mathcal{R}}$. Associated to a positive root 
$\widehat{\alpha} \in \widehat{\Pi }$ is a root algebra ${\widetilde{\mathfrak{g}}}^{(\widehat{\alpha })}$, which is isomorphic to a $3$-dimensional Heisenberg algebra ${\mathrm{h}}_3$. In \S 2.3 we associate to every extremal root chain ${\mathcal{S}}^{\widehat{\beta }}_{\widehat{\alpha }}$ the space 
${\widehat{\mathfrak{g}}}^{\widehat{\beta }}_{\widehat{\alpha }}$,  
formed by taking the direct sum of root spaces corresponding to elements of 
${\mathcal{S}}^{\widehat{\beta }}_{\widehat{\alpha }}$. In claim 2.2.1 we show that the adjoint representation of 
${\widetilde{\mathfrak{g}}}^{(\widehat{\alpha })}$ on $\widetilde{\mathfrak{n}}$ is completely reducible and classify its irreducible summands. The goal of \S 2.2 is to verify that $\Phi $ satisfies axioms 4 and 5 of a system of roots. Lemma 2.3.1 shows that $\Phi $ satisfies axiom $5$ of a system of roots. 
The main result of this section is claim 2.3.2, which shows that for every fixed $\widehat{\alpha } \in \Phi $ the function $K_{\widehat{\alpha }}: \widehat{\mathcal{R}} \rightarrow \Z : \widehat{\beta }\mapsto \langle \widehat{\beta} , \widehat{\alpha } \rangle $ is linear, that is, $\Phi $ satisfies axiom $4$ of a system of roots. Our main tool in proving claim 2.3.2 is claim 2.3.4, which states that if for $i=1,2$ 
the extremal root chain ${\mathcal{S}}^{{\widehat{\beta }}_i}_{\widehat{\alpha }}$ 
has an integer pair $(q_i,p_i)$ and if ${\mathcal{S}}^{{\widehat{\beta }}_1 +  
{\widehat{\beta }}_2}_{\widehat{\alpha }}$ is also an extremal root chain with integer 
pair $(r,s)$, then the Killing integer $r-s$ of 
${\mathcal{S}}^{{\widehat{\beta }}_1+ {\widehat{\beta }}_2}_{\widehat{\alpha }}$ 
is equal to $(q_1+q_2)-(p_1+p_2)$.  

\subsection{Roots and root spaces}
\label{sec2subsec1}

Let $\widetilde{\mathfrak{g}} = \mathfrak{g} \oplus \widetilde{\mathfrak{n}}$ be a very special 
sandwich algebra with nilpotent radical $\widetilde{\mathfrak{n}}$ and $\mathfrak{g}$ a complex simple Lie algebra with Cartan subalgebra $\mathfrak{h}$. \medskip

\noindent \textbf{Lemma 2.1.1} The center $Z$ of the nilpotent radical $\widetilde{\mathfrak{n}}$ is an 
${\ad }_{\mathfrak{h}}$-invariant subspace of $\widetilde{\mathfrak{n}}$. \medskip 

\noindent \textbf{Proof.} To see this let $H\in \mathfrak{h}$, $Y_1 \in \widetilde{\mathfrak{n}}$, 
and $Z_1 \in Z$. Then 
\begin{equation}
[ [H,Z_1], Y_1] = [ [ H,Y_1], Z_1] + [H, [Z_1,Y_1]] =0 
\label{eq-s2onenw}
\end{equation}
holds, since $[H, Y_1] \in \widetilde{\mathfrak{n}}$ and $Z_1 \in Z$ implies that $[[H,Y_1], Z_1] =0$ and 
$[H, [Z_1,Y_1]] =0$, since $Z_1\in Z$ and $Y_1  \in \widetilde{\mathfrak{n}}$ implies $[Z_1,Y_1] =0$. 
Now $[H, Z_1] \in \widetilde{\mathfrak{n}}$ 
because ${\ad }_H$ maps $\widetilde{\mathfrak{n}}$ into itself and $Z_1 \in Z \subseteq 
\widetilde{\mathfrak{n}}$. From equation (\ref{eq-s2onenw}) it follows that $[H,Z_1]$ 
commutes with every element of $\widetilde{\mathfrak{n}}$. Thus $[H, Z_1] \in Z$. So ${\ad }_H$ maps $Z$ into itself for every $H \in \mathfrak{h}$, as desired. \hfill $\square $ \medskip 

Let $\mathcal{Z} = \{ \zeta \in \mathfrak{R} \} $ and set $\mathcal{Y} = \widehat{\mathcal{R}} = 
\mathfrak{R} \setminus \mathcal{Z}$. Since ${\ad }_{\mathfrak{h}}|\widetilde{\mathfrak{n}}$ is 
a maximal torus of $\gl (\widetilde{\mathfrak{n}}, \C)$, its set of root vectors 
$\{ X_{\tau } \in \widetilde{\eta } \setrule \, \tau \in \mathfrak{R} \}$ is a basis of 
$\widetilde{\mathfrak{n}}$ consisting of eigenvectors of 
${\ad }_H|\widetilde{\mathfrak{n}}$ for every $H \in \mathfrak{h}$. Since $Z$ is ${\ad }_{\mathfrak{h}}$ 
invariant, because $\widetilde{\mathfrak{n}}$ is an ideal of $\widetilde{\mathfrak{g}}$ and 
$Z \subseteq \widetilde{\mathfrak{n}}$, there is an ${\ad }_{\mathfrak{h}}$ invariant subspace $Y$ of $\widetilde{\mathfrak{n}}$, 
which is complementary to $Z$. So $Y = \sum_{\widehat{\alpha } \in \mathcal{Y}} \oplus 
{\widetilde{\mathfrak{n}}}_{\widehat{\alpha }}$, where 
\begin{displaymath}
{\widetilde{\mathfrak{n}}}_{\widehat{\alpha }} = \{ X \in Y \setrule \, 
[H, X] = \widehat{\alpha }(H) X \, \, \mbox{for all $H \in \mathfrak{h}$} \} . 
\end{displaymath}
${\widehat{\mathfrak{n}}}_{\widehat{\alpha }}$ is the \emph{root space} in $Y$ corresponding to the root 
$\widehat{\alpha } \in \widehat{\mathcal{R}}$. 
Since ${\ad }_{\mathfrak{h}}|\widetilde{\mathfrak{n}}$ is a maximal torus, 
${\dim }_{\C}\, {\widehat{\mathfrak{n}}}_{\widehat{\alpha }}  =1$. So 
${\widehat{\mathfrak{n}}}_{\widehat{\alpha }}$ is spanned by a nonzero vector $X_{\widehat{\alpha}}$, 
which is called a \emph{root vector} in $\widetilde{\mathfrak{n}}$ associated to the root $\widehat{\alpha }$. \medskip   

\noindent \textbf{Claim 2.1.2} For every $\widehat{\alpha } \in \widehat{\mathcal{R}}$ we have 
\begin{equation}
{\widetilde{\mathfrak{n}}}_{\widehat{\alpha }} = {\underline{\mathfrak{g}}}_{\underline{\alpha }}, 
\label{eq-s2zeronw*}
\end{equation}
where ${\underline{\mathfrak{g}}}_{\underline{\alpha }}$ is the root space in $\underline{\mathfrak{g}}$ 
with Cartan subalgebra $\underline{\mathfrak{h}}$ corresponding to the root $\underline{\alpha} \in \underline{\mathcal{R}}$ such that 
$\underline{\alpha }|\mathfrak{h} = \widehat{\alpha }$ and $\underline{\alpha }(\underline{\widetilde{H}}) 
< 0$. \medskip 

\noindent \textbf{Proof.} Since $\widetilde{\mathfrak{n}}$ is an ideal in $\underline{\mathfrak{g}}$, 
$\widetilde{\mathfrak{n}}$ is an ${\ad}_{\underline{\mathfrak{h}}}$-invariant subalgebra of 
$\underline{\mathfrak{g}}$. 
Let $X_{\widehat{\alpha}}$ be a root vector in $\widetilde{\mathfrak{n}}$ corresponding to the 
root $\widehat{\alpha } \in \widehat{\mathcal{R}}$. Because $\widetilde{\mathfrak{n}} = 
\sum_{\widehat{\beta } \in \mathfrak{R}} \oplus {\widetilde{\mathfrak{n}}}_{\widehat{\beta }}$, we get 
\begin{align}
[\underline{H}, X_{\widehat{\alpha }}] & = 
\sum_{\widehat{\beta } \in \mathfrak{R}} c_{\widehat{\beta }}(\underline{H}) X_{\widehat{\beta }} 
\in \widetilde{\mathfrak{n}} , 
\label{eq-s2twonw*}
\end{align}
for every $\underline{H} \in \underline{\mathfrak{h}}$. Here $c_{\widehat{\beta }} \in 
{\underline{\mathfrak{h}}}^{\ast }$ for every $\widehat{\beta } \in \mathfrak{R}$. If $[\underline{H}, X_{\widehat{\alpha }}] =0$ for every $\underline{H} \in \underline{\mathfrak{h}}$, then 
${\spann}_{\C}\{ \underline{H}, \underline{H}\in \mathfrak{h}; X_{\widehat{\alpha}} \} $ is 
an abelian subalgebra of $\underline{\mathfrak{g}}$, which properly contains $\underline{\mathfrak{h}}$. 
This contradicts the fact that $\underline{\mathfrak{h}}$ is a Cartan subalgebra of 
$\underline{\mathfrak{g}}$. Hence $[\underline{H}, X_{\widehat{\alpha}}] \ne 0$. So 
$c_{\widehat{\beta}} \ne 0$ for some $\widehat{\beta } \in \mathfrak{R}$. Since $\mathfrak{h} 
\subseteq \underline{\mathfrak{h}}$, equation (\ref{eq-s2twonw*}) yields 
\begin{align}
[H, X_{\widehat{\alpha}}] & = \sum_{\widehat{\beta } \in \mathfrak{R}} c_{\widehat{\beta }}(H) X_{\widehat{\beta }} = \widehat{\alpha }(H) X_{\widehat{\alpha }}, 
\label{eq-s2threenw*}
\end{align}
for every $H \in \mathfrak{h}$. Thus $c_{\widehat{\beta }}|\mathfrak{h} =0$, if $\widehat{\beta } \ne 
\widehat{\alpha }$ and $c_{\widehat{\alpha}}|\mathfrak{h} = \widehat{\alpha }$. Since 
$\underline{\mathfrak{h}} = \mathfrak{h} \oplus {\spann}_{\C}\{ \underline{\widetilde{H}} \}$, it 
follows that $c_{\widehat{\beta }} = c_{\widehat{\beta }}|{\spann}_{\C}\{ \underline{\widetilde{H}} \}$ for every 
$\widehat{\beta} \ne \widehat{\alpha }$. From 
(\ref{eq-s2twonw*}) we get 
\begin{align}
[\underline{\widetilde{H}}, X_{\widehat{\alpha }}] & = \sum_{\widehat{\beta } \in \mathfrak{R}} c_{\widehat{\beta }}(\underline{\widetilde{H}}) X_{\widehat{\beta }} \in \widetilde{\mathfrak{n}}. 
\label{eq-s2fournw*}
\end{align}
For every $H \in \mathfrak{h}$ we have 
\begin{align}
[H, [\underline{\widetilde{H}}, X_{\widehat{\alpha}}]] & = 
[ [ H, \underline{\widetilde{H}}], X_{\widehat{\alpha}} ] + [\underline{\widetilde{H}}, [ H, X_{\widehat{\alpha}}]] 
= \widehat{\alpha }(H) [\underline{\widetilde{H}}, X_{\widehat{\alpha }}], 
\label{eq-s2fivenw*}
\end{align}
since $[H, \underline{\widetilde{H}}] =0$, because $H$, $\underline{\widetilde{H}} \in \underline{\mathfrak{h}}$ and $\underline{\mathfrak{h}}$ is an abelian subalgebra of $\underline{\mathfrak{g}}$. If 
$[\underline{\widetilde{H}}, X_{\widehat{\alpha }}] = 0$, then by equation 
(\ref{eq-s2fournw*}) $c_{\widehat{\beta }}(\underline{\widetilde{H}}) =0$ for every 
$\widehat{\beta } \in \mathfrak{R}$. Thus $c_{\widehat{\beta }} =0$ for every 
$\widehat{\beta } \in \mathfrak{R}$ such that $\widehat{\beta } \ne \widehat{\alpha }$ and $c_{\widehat{\alpha }}(\underline{\widetilde{H}}) =0$. From equation (\ref{eq-s2twonw*}) we obtain 
\begin{equation}
[\underline{H}, X_{\widehat{\alpha }}] = c_{\widehat{\alpha}}(\underline{H}) \, X_{\widehat{\alpha}} \quad 
\mbox{for every $\underline{H} \in \underline{\mathfrak{h}}$.}
\label{eq-s2sixnw*} 
\end{equation}
So $c_{\widehat{\alpha}} \in \underline{\mathcal{R}}$. Since $c_{\widehat{\alpha }}|\mathfrak{h} = 
\widehat{\alpha}$, the fact that $c_{\widehat{\alpha}}(\underline{\widetilde{H}}) =0$ contradicts 
the hypothesis that $\widehat{\alpha } \in \widehat{\mathcal{R}}$. Thus $[ \underline{\widetilde{H}}, 
X_{\widehat{\alpha }}] \ne 0$. From equation (\ref{eq-s2fivenw*}) 
it follows that $[\underline{\widetilde{H}}, X_{\widehat{\alpha }} ]$ is a root vector in $\widetilde{\mathfrak{n}}$ 
corresponding to the root $\widehat{\alpha } \in \mathfrak{R}$, that is, 
\begin{displaymath}
[\underline{\widetilde{H}}, X_{\widehat{\alpha}}] = 
d\,  X_{\widehat{\alpha }}, 
\end{displaymath}
for some $d \in \R$. The number $d$ is nonzero, since $X_{\widehat{\alpha }}$ is a basis of 
${ \widetilde{\mathfrak{n}} }_{\widehat{\alpha }}$, and 
$[ \underline{\widetilde{H}}, X_{\widehat{\alpha }}] \ne 0$. From equation 
(\ref{eq-s2twonw*}) we obtain equation (\ref{eq-s2sixnw*}), which implies that 
$c_{\widehat{\alpha }}(\underline{\widetilde{H}}) = d \ne 0$ and $c_{\widehat{\alpha}} \in 
\underline{\mathcal{R}}$, whose root space ${\underline{\mathfrak{g}}}_{c_{\widehat{\alpha}}}$ is spanned 
by $X_{\widehat{\alpha}}$. Now $c_{\widehat{\alpha }}|\mathfrak{h} = \widehat{\alpha}$. 
So $c_{\widehat{\alpha }}(\underline{\widetilde{H}}) < 0$, since $c_{\widehat{\alpha }} \in 
\widehat{\mathcal{R}}$. Thus $c_{\widehat{\alpha }}$ is the 
root in $\underline{\mathcal{R}}$ which we previously called $\underline{\alpha }$. Therefore 
equation (\ref{eq-s2zeronw*}) holds. \hfill $\square $ \medskip

Define an addition operation $+$ on $\widehat{\mathcal{R}}$ by saying that if 
$\widehat{\alpha }, \widehat{\beta } \in \widehat{\mathcal{R}}$, then 
$\widehat{\alpha } + \widehat{\beta } \in \widehat{\mathcal{R}}$ if there are 
$\underline{\alpha }$, $\underline{\beta } \in \underline{\mathcal{R}}$ such that 
$\underline{\alpha } + \underline{\beta } \in \underline{\mathcal{R}}$ and $\widehat{\alpha } = 
\underline{\alpha } | \mathfrak{h}$ and  $\widehat{\beta } = \underline{\beta } | \mathfrak{h}$. In 
particular, addition of $\widehat{\alpha }$ and $\widehat{\beta }$ in $\widehat{\mathcal{R}}$ is 
just the addition of the values of the covectors $\widehat{\alpha }$ and $\widehat{\beta }$ in 
$\mathfrak{h}$. Clearly the operation of addition on $\widehat{\mathcal{R}}$ is 
commutative and associative. \medskip 

Next we prove \medskip

\noindent \textbf{Claim 2.1.3} For every $\widehat{\alpha }$, $\widehat{\beta } \in \widehat{\mathcal{R}}$
\begin{equation}
[ {\widetilde{\mathfrak{n}}}_{\widehat{\alpha}}, {\widetilde{\mathfrak{n}}}_{\widehat{\beta}} ] = 
\left\{ \begin{array}{rl} 
{\widetilde{\mathfrak{n}}}_{\widehat{\alpha} + \widehat{\beta }} , & 
\mbox{if $\widehat{\alpha } +\widehat{\beta } \in \widehat{\mathcal{R}}$} \\
\{ 0 \} , & \mbox{otherwise.} 
\end{array} \right. 
\label{eq-s2sevennw*}
\end{equation}

\noindent \textbf{Proof.} If $\widehat{\alpha } +\widehat{\beta } \in \widehat{\mathcal{R}}$, 
then $[ {\widetilde{\mathfrak{n}}}_{\widehat{\alpha}}, {\widetilde{\mathfrak{n}}}_{\widehat{\beta}} ] 
\subseteq {\widetilde{\mathfrak{n}}}_{\widehat{\alpha} + \widehat{\beta }}$, since for every $H \in 
\mathfrak{h}$
\begin{align*}
[ H, [X_{\widehat{\alpha}}, X_{\widehat{\beta }} ] ] & = 
[ [H, X_{\widehat{\alpha}}], X_{\widehat{\beta }}] + [X_{\widehat{\alpha}}, [ H, X_{\widehat{\beta }}] ] \\
& = \widehat{\alpha}(H)[X_{\widehat{\alpha}}, X_{\widehat{\beta }} ] + 
\widehat{\beta }(H) [X_{\widehat{\alpha}}, X_{\widehat{\beta }} ] \\
& = (\widehat{\alpha } + \widehat{\beta })(H) [X_{\widehat{\alpha}}, X_{\widehat{\beta }} ] ,
\end{align*}
which implies $[X_{\widehat{\alpha}}, X_{\widehat{\beta }} ] \in 
{\widetilde{\mathfrak{n}}}_{\widehat{\alpha} + \widehat{\beta }}$. To show that 
$[ {\widetilde{\mathfrak{n}}}_{\widehat{\alpha}}, {\widetilde{\mathfrak{n}}}_{\widehat{\beta}} ] = 
{\widetilde{\mathfrak{n}}}_{\widehat{\alpha} + \widehat{\beta }}$, we must verify that the 
vector $[X_{\widehat{\alpha}}, X_{\widehat{\beta }} ]$ is nonzero and thus forms a basis of 
the root space ${\widetilde{\mathfrak{n}}}_{\widehat{\alpha}+ \widehat{\beta }}$. Since $\widehat{\alpha }$, 
$\widehat{\beta } \in \widehat{\mathcal{R}}$, there are $\underline{\alpha }$, $\underline{\beta } 
\in \underline{\mathcal{R}}$ such that $\underline{\alpha }|\mathfrak{h} = \widehat{\alpha }$, 
$\underline{\beta }|\mathfrak{h} = \widehat{\beta }$ and $\underline{\alpha }(\underline{\widetilde{H}}) < 0$, 
$\underline{\beta }(\underline{\widetilde{H}}) < 0$, $\underline{\alpha} + \underline{\beta } \in 
\underline{\mathcal{R}}$, and $\widehat{\alpha } + \widehat{\beta } = 
(\underline{\alpha } + \underline{\beta })|\mathfrak{h}$. By claim 2.1.2 we have 
${\widetilde{\mathfrak{n}}}_{\widehat{\alpha}} = {\underline{\mathfrak{g}}}_{\underline{\alpha }}$ and  
${\widetilde{\mathfrak{n}}}_{\widehat{\beta}} = {\underline{\mathfrak{g}}}_{\underline{\beta }}$. Thus 
$X_{\widehat{\alpha }}$ and $X_{\widehat{\beta }}$ are root vectors in the root spaces 
${\underline{\mathfrak{g}}}_{\underline{\alpha }}$ and ${\underline{\mathfrak{g}}}_{\underline{\beta }}$, 
respectively. Since $\underline{\mathfrak{g}}$ is a simple Lie algebra, 
$[X_{\widehat{\alpha }}, X_{\widehat{\beta }}] \ne 0$, see Humphreys \cite[p.39]{humphreys}. \medskip 

If $\widehat{\alpha } + \widehat{\beta } \notin \widehat{\mathcal{R}}$, then 
$[X_{\widehat{\alpha }}, X_{\widehat{\beta }}] = 0$, because $\underline{\alpha } + \underline{\beta } 
\notin \underline{\mathcal{R}}$, which implies ${\underline{\mathfrak{g}}}_{\underline{\alpha } + \underline{\beta }} = \{ 0 \}$. Thus ${\widetilde{\mathfrak{n}}}_{\widehat{\alpha } + \widehat{\beta }} = \{ 0 \}$. 
\hfill $\square $ \medskip 

We now show that \medskip 

\noindent \textbf{Lemma 2.1.4} Given $\widehat{\alpha } \in \mathcal{Y}= \widehat{\mathcal{R}}$ there is 
$\widehat{\beta } \in \mathcal{Y}$ satisfying $\widehat{\alpha } + \widehat{\beta } = \zeta $. \medskip 

\noindent \textbf{Proof.} Suppose that there is a $\widehat{\gamma } \in \widehat{\mathcal{R}}$ 
such that $[X_{\widehat{\gamma }}, X_{\widehat{\beta }} ] =0$ for every $\widehat{\beta } \in \mathcal{Y}$. 
Since $Z$ is the center of $\widetilde{\mathfrak{n}}$, it follows that $[X_{\zeta }, X_{\widehat{\gamma }}] = 0$. 
Thus $X_{\widehat{\gamma }}$ commutes with every element of $\widetilde{\mathfrak{n}} = 
{\spann }_{\C } \{ X_{\zeta }, X_{\widehat{\beta }}, \, \widehat{\beta } \in \mathcal{Y} \}$. So $X_{\widehat{\gamma }} \in Z$, which gives $\widehat{\gamma } \in \mathcal{Z}$. But by 
hypothesis $\widehat{\gamma } \in \mathcal{Y}$. Since $\mathcal{Z} \cap \mathcal{Y} = \{ 0 \}$, 
we obtain a contradiction. Thus for each $\widehat{\alpha } \in \mathcal{Y}$ there is 
a $\widehat{\beta } \in \mathcal{Y}$ such that $[X_{\widehat{\alpha }}, X_{\widehat{\beta }}] \ne 0$. 
Since $\widetilde{\mathfrak{n}}$ is a sandwich, $[Y,Y] \subseteq Z$. Because 
$\{X_{\zeta }\} $ is a basis of $Z$, there is a $d \in \C \setminus \{ 0 \}$ such that 
$[X_{\widehat{\alpha }}, X_{\widehat{\beta }} ] = d \, X_{\zeta}$. Rescaling the root vectors $X_{\widehat{\alpha}}$ and 
$X_{\widehat{\beta }}$ we may assume that $[X_{\widehat{\alpha }}, X_{\widehat{\beta }} ] =  X_{\zeta}$. For 
every $H \in \mathfrak{h}$
\begin{align*}
\zeta (H) X_{\zeta } & = [H, X_{\zeta }] = [H, [X_{\widehat{\alpha }}, X_{\widehat{\beta }} ] ] 
= [ [ H, X_{\widehat{\alpha}} ], X_{\widehat{\beta }} ] + [X_{\widehat{\beta }}, [H, X_{\widehat{\alpha }} ] ] 
\\
& = \big( \widehat{\alpha }(H) + \widehat{\beta }(H)\big) [X_{\widehat{\alpha }}, X_{\widehat{\beta }} ] 
= \big( \widehat{\alpha }(H) + \widehat{\beta }(H)\big) X_{\zeta}. 
\end{align*}
Since $X_{\zeta } \ne 0$, we get $\widehat{\alpha }(H) + \widehat{\beta }(H) = \zeta (H)$ for every 
$H \in \mathfrak{h}$, that is, $\widehat{\alpha } + \widehat{\beta } = \zeta $. \hfill $\square $ \medskip 

\noindent \textbf{Corollary 2.1.4A} The covector $\widehat{\beta } \in \mathcal{Y}$ satisfying 
the hypothesis of lemma 2.1.4 is unique. \medskip 

\noindent \textbf{Proof.} Suppose that in addition to $\widehat{\beta } \in \widehat{R}$ such that 
$\widehat{\alpha } + \widehat{\beta } = \zeta $, there is a ${\widehat{\beta }}\, ' \in \widehat{R}$ such that 
$\widehat{\alpha } + {\widehat{\beta }}\, ' = \zeta $. Then 
\begin{displaymath}
{\widehat{\beta }}\, ' + \zeta = {\widehat{\beta }}\, ' + (\widehat{\alpha } + \widehat{\beta }) =
(\widehat{\alpha } + {\widehat{\beta }}\, ' ) + \widehat{\beta } = \zeta + \beta . 
\end{displaymath}
So for every $H \in \mathfrak{h}$ we get 
${\widehat{\beta }}\,'(H) + \zeta (H) = \zeta (H) + \widehat{\beta } (H)$, 
which implies ${\widehat{\beta }}\, '(H) = \widehat{\beta}(H)$ for every $H \in \mathfrak{h}$, since 
$\zeta (H) = 0$ for all $H \in \mathfrak{h}$. Therefore ${\widehat{\beta }}\, ' = \widehat{\beta}$. 
\hfill $\square $ \medskip 

Since $\zeta $ is the zero covector in ${\mathfrak{h}}^{\ast}$, the root $\widehat{\beta }$ given 
by lemma 2.1.4 is the negative of the given root $\widehat{\alpha } \in \widehat{\mathcal{R}}$, 
which we denote by $-\widehat{\alpha}$. Now 
$-\widehat{\alpha } + \big( - ( - \widehat {\alpha }) \big) = \zeta = -\widehat{\alpha} + \widehat{\alpha }$. 
So using corollary 2.1.4A we find that $\widehat{\alpha } = - ( - \widehat {\alpha })$. Thus 
$-\widehat{\mathcal{R}} = \widehat{\mathcal{R}}$. Consequently, axioms 1 and 2 of a system of 
roots holds for $V = Y$ with $\Phi = \widehat{\mathcal{R}}$. Therefore by lemma 1.2.1 
axiom 3 of a system of roots holds. \medskip   

Let ${\underline{\Pi }}^{+}$ be the set of simple positive roots in $\underline{\mathcal{R}} \subseteq 
{\underline{\mathfrak{h}}}^{\ast }$ for the simple Lie algebra $\underline{\mathfrak{g}}$. Let 
${\underline{\widehat{\Pi}}}^{+}$ be the set of $\widehat{\alpha } \in \widehat{\mathcal{R}}$ such that 
there is an $\underline{\alpha } \in {\underline{\Pi }}^{+}$ with $\underline{\alpha}|\mathfrak{h} = 
\widehat{\alpha }$. \medskip 

\noindent \textbf{Lemma 2.1.5} Every $\widehat{\alpha } \in {\underline{\widehat{\Pi }}}^{+}$ is simple. \medskip 

\noindent \textbf{Proof.} Given $\widehat{\alpha } \in {\underline{\widehat{\Pi }}}^{+}$, suppose that 
$\widehat{\alpha } = \widehat{\beta } + \widehat{\gamma }$ for 
some $\widehat{\beta }$, $\widehat{\gamma } \in \widehat{\mathcal{R}}$. Then there are 
$\underline{\alpha} \in {\underline{\Pi }}^{+}$ and $\underline{\beta }$, $\underline{\gamma} \in 
\underline{\mathcal{R}}$ such that $\widehat{\alpha } = \underline{\alpha }|\mathfrak{h}$, 
$\widehat{\beta } = \underline{\beta }|\mathfrak{h}$, and $\widehat{\gamma } = 
\underline{\gamma }|\mathfrak{h}$ with $\underline{\alpha }|\mathfrak{h} = \underline{\beta }|\mathfrak{h} 
+\underline{\gamma}|\mathfrak{h}$. Since $\widehat{\alpha }$, $\widehat{\beta }$, and 
$\widehat{\gamma } \in \widehat{\mathcal{R}} \subseteq {\mathfrak{h}}^{\ast }$ and 
$\widetilde{\mathfrak{g}}$ is a very special sandwich algebra, it follows that 
$\underline{\alpha}(\underline{\widetilde{H}}) = 
\underline{\beta }(\underline{\widetilde{H}}) = \underline{\gamma }(\underline{\widetilde{H}}) =0$. From   
$\underline{\mathfrak{h}} = \mathfrak{h} \oplus {\spann }_{\C}\{ \underline{\widetilde{H}} \}$ we obtain 
$\underline{\alpha } = \underline{\beta } + \underline{\gamma }$. This contradicts the fact that 
$\underline{\alpha } \in {\underline{\Pi }}^{+}$. Thus $\widehat{\alpha }$ is a simple root in 
$\widehat{\mathcal{R}}$. \hfill $\square $ \medskip 

Let $\underline{\Pi } = {\spann}_{{\R }_{>0}}\{ \underline{\alpha } \in {\underline{\Pi }}^{+} \} $. Then 
$\underline{\Pi }$ is the set of positive roots in $\underline{\mathcal{R}}$. Let $\widehat{\Pi }$ be 
the set of $\widehat{\alpha } \in \widehat{\mathcal{R}}$ such that there is an $\underline{\alpha } \in 
\underline{\Pi}$ with $\underline{\alpha }|\mathfrak{h} = \widehat{\alpha }$. Then 
$\widehat{\Pi } = {\spann }_{{\R }_{>0}}\{ \widehat{\alpha } \in {\widehat{\Pi}}^{+} \}$ is the 
set of \emph{positive} roots in $\widehat{\mathcal{R}}$. \medskip 

\noindent \textbf{Claim 2.1.6} $\widehat{\mathcal{R}} = \mathcal{Y} = \widehat{\Pi } \oplus (-\widehat{\Pi })$. \medskip 

\noindent \textbf{Proof.} Suppose that $\widehat{\alpha } \in \mathcal{Y}$. Then there 
is an $\underline{\alpha } \in \underline{\mathcal{R}} \subseteq {\underline{\mathfrak{h}}}^{\ast}$ 
such that $\widehat{\alpha } = \underline{\alpha }|\mathfrak{h}$. Since $\underline{\mathfrak{g}}$ is 
simple, we have $\underline{\mathfrak{g}} = \underline{\Pi } \oplus (-\underline{\Pi })$. Thus there 
are $\underline{\beta } \in \underline{\Pi}$ and $\underline{\gamma } \in -(\underline{\Pi })$ such that  
$\underline{\alpha } = \underline{\beta } + \underline{\gamma}$. Set $\widehat{\beta } = 
\underline{\beta }|\mathfrak{h}$ and $\widehat{\gamma } = \underline{\gamma }|\mathfrak{h}$. 
Then $\widehat{\beta } \in \widehat{\Pi}$, $\widehat{\gamma } \in -(\widehat{\Pi})$, and $\widehat{\alpha } = \widehat{\beta } + \widehat{\gamma }$. So $\mathcal{Y} \subseteq \widehat{\Pi } 
\oplus (-\widehat{\Pi })$. The reverse inclusion $\widehat{\Pi} \oplus (-\widehat{\Pi }) \subseteq 
\mathcal{Y}$ is immediate. Hence $\mathcal{Y} = \widehat{\Pi } \oplus (-\widehat{\Pi })$. \hfill $\square $ 

\subsection{Root subalgebras}
\label{sec2subsec2}

For each root $\widehat{\alpha } \in \widehat{\Pi}$, the \emph{root subalgebra} 
${\widetilde{\mathfrak{n}}}^{(\widehat{\alpha})} = {\spann }_{\C} \{ X_{\widehat{\alpha }}, \,  
X_{-\widehat{\alpha}}, \, H_{\zeta} = X_{\zeta } \} $ is the Heisenberg algebra ${\mathrm{h}}_3$ with bracket relations
\begin{displaymath}
[H_{\zeta }, X_{\widehat{\alpha }}] = 
[H_{\zeta}, X_{-\widehat{\alpha} }] \, = \, 0 \, \, \, \mathrm{and}\,  \, \,   
[X_{\widehat{\alpha }}, X_{-\widehat{\alpha }}] =  H_{\zeta }. 
\end{displaymath}
First we prove \medskip

\noindent \textbf{Claim 2.2.1} Every finite dimensional representation of 
${\mathrm{h}}_3$ by nilpotent \linebreak 
linear maps of a finite dimensional complex vector space $V$ into 
itself is completely reducible. \medskip 

\noindent \textbf{Proof.} Let $\{ X,Y, H \} $ span ${\mathrm{h}}_3$ with bracket relations 
\begin{displaymath}
[H, X] = [H,Y] = 0, \, \, \, \mathrm{and} \, \, \, [X,Y] = H, 
\end{displaymath} 
Let $\rho  : {\mathrm{h}}_3 \rightarrow \gl (V, \C )$ be a representation of Lie algebras 
such that $\rho (Y)$, $\rho (X)$, and $\rho (H)$ are complex nilpotent linear maps of $V$ into itself. 
For $v \in \ker \rho (Y)$ let $\mathfrak{f} = \{ w, \rho (Y)w, \ldots , ({\rho }(Y))^n w =v \}$ be the longest Jordan chain of $\rho (Y)$ in $V$, which ends at $v$. Suppose that $n \ge 2$. Then the matrix of 
$\rho (Y)$ with respect to the basis $\mathfrak{f}$ of the vector space $U = 
{\spann }_{\C} \{ w, \rho (Y)w, \ldots , ({\rho }(Y))^n w = 
v \} $ is the lower $(n+1) \times (n+1)$ Jordan block. With respect to the basis $\mathfrak{f}$ the 
$(n+1) \times (n+1)$ matrices of $\rho (Y)$, $\rho (X)$, and $\rho (H)$ are  
\begin{displaymath} 
\mbox{{\tiny $\begin{pmatrix}
0 & &  & & & & \\ 1 & 0 & & && & \\ & & & & & & \\ 
 & \ddots  & & & & & & \\ & & & & \ddots & & \\
 & & & & \ddots  & 0&   \\ & & & &  & 1 &  0 \end{pmatrix} $}}, \, 
\mbox{{\tiny $\begin{pmatrix} 
0 & & & & & \\ 
\vdots & \ddots & & & & \\
\vdots & & \ddots  & &  & \\ 
0 & & & \ddots &  & \\
-1 & 0 & \cdots & \cdots & 0 & \\ 
0 & 1 & 0 &\cdots & 0 & 0 \end{pmatrix} $}} , \, 
\mbox{{\tiny $\left( \begin{array}{ccclr} 0 & & & & \\
0 & \ddots & & & \\ 
\vdots & \ddots & \ddots  & &  \\ 
0 & \cdots & \ddots & 0 & \\ 
2 & 0 &  \cdots & 0 & 0 
\end{array} \right)  $,}}
\end{displaymath}
respectively. When $n=0$ let $\rho (Y)$, $\rho (X)$, 
and $\rho (H)$ be the zero matrix; while when $n=1$ let $\rho (Y) =${\tiny 
$\begin{pmatrix} 0 & 0 \\ 1 & 0 \end{pmatrix}$}, and let $\rho (X)$ and 
$\rho (H)$ be the $2\times 2$ zero matrix. For every $n \ge 0$ it is straightforward 
to check that the following bracket relations hold
\begin{displaymath}
[\rho (H), \rho (X) ] = [\rho (H), \rho (Y) ] = 0 \, \, \, 
\mathrm{and} \, \, \, [\rho (X), \rho (Y)] = \rho (H).
\end{displaymath} 
This verifies that the complex linear maps $\rho (Y)$, $\rho (X)$, and $\rho (H)$ of $U$ into itself 
form a Lie subalgebra of $\gl (V, \C)$, which is isomorphic to ${\mathrm{h}}_3$. This representation 
of ${\mathrm{h}}_3$ is irreducible. Suppose not. 
Then there is a $\rho ({\mathrm{h}}_3)$-invariant subspace $W$ of $U$ such that 
$U = W\oplus {\spann}_{\C}\{ v \} $. But $(\rho (Y))^n W \subseteq \spann \{ v \} $, 
which contradicts the $\rho ({\mathrm{h}}_3)$-invariance of $W$. Therefore on $U$ 
the representation $\rho  $ is irreducible. Repeating the above construction for each vector in a basis 
of $\ker \rho (Y)$ determines the Jordan normal form of $\rho (Y)$ and thus decomposes the representation 
$\rho $ into a sum of finite dimensional irreducible representations 
of ${\mathrm{h}}_3$. The summands in this representation are unique up to a reordering, 
because the Jordan blocks of $\rho (Y)$ are unique up to a reordering. \hfill $\square $ \medskip 

Let $\widehat{\alpha } \in \widehat{\Pi}$. If $\widehat{\alpha } = \zeta $, then for every 
$\widehat{\beta } \in \mathfrak{R}$ we have $[X_{\widehat{\alpha}}, X_{\widehat{\beta }}] \in 
[Z ,\widetilde{\mathfrak{n}}] = \{ 0 \}$. Thus ${\ad }_{X_{\zeta}}$ is a nilpotent complex linear 
map of $\widetilde{\mathfrak{n}}$ into itself. Let $\widehat{\alpha } \ne \zeta$. 
From claim 2.1.3 it follows that ${\ad }_{X_{\widehat{\alpha}}}$ and ${\ad}_{X_{-\widehat{\alpha}}}$ 
are complex nilpotent linear maps of $\widetilde{\mathfrak{n}} = Z \oplus Y = 
{\spann}_{\C}\{ X_{\zeta } \} \oplus \sum_{\widehat{\beta } \in \widehat{\mathcal{R}}} \oplus 
{\widetilde{\mathfrak{n}}}_{\widehat{\beta }}$. In more detail, for every $\widehat{\beta } \in \widehat{\mathcal{R}}$ we have 
\begin{displaymath}
{\ad }_{X_{\widehat{\alpha }}}X_{\widehat{\beta }} = \left\{ \begin{array}{rl}
X_{\widehat{\beta} + \widehat{\alpha}}, & \mbox{if $\widehat{\beta } +\widehat{\alpha } \in 
\widehat{\mathcal{R}}$} \\
0, & \mbox{otherwise.} 
\end{array} \right. 
\end{displaymath}
So ${\ad}^n_{X_{\widehat{\alpha}}}X_{\widehat{\beta}} = \{ 0 \}$ as soon as $\widehat{\beta } + 
n \widehat{\alpha } \notin \widehat{\mathcal{R}}$. Such an integer $n$ exists and is finite because 
$\widehat{\mathcal{R}}$ is a finite set. Thus ${\ad }_{X_{\widehat{\alpha}}}$ is a complex 
nilpotent linear mapping of $\widetilde{\mathfrak{n}}$ into itself. A similar argument shows that 
${\ad }_{X_{-\widehat{\alpha}}}$ is also a nilpotent linear mapping of $\widetilde{\mathfrak{n}}$ 
into itself. This proves \medskip 

\noindent \textbf{Lemma 2.2.2} For each $\widehat{\alpha } \in \widehat{\Pi}$ the mapping 
$\ad : {\widetilde{\mathfrak{n}}}^{(\widehat{\alpha })} \rightarrow \gl (\widetilde{\mathfrak{n}}, \C)$, 
which sends $X_{\widehat{\alpha}}$, $X_{-\widehat{\alpha}}$, and $H_{\zeta } = X_{\zeta }$ to 
the nilpotent linear maps ${\ad }_{X_{\widehat{\alpha}}}$, ${\ad }_{-\widehat{\alpha}}$, and 
${\ad }_{X_{\zeta}}$ of $\widetilde{\mathfrak{n}}$ into itself, respectively, is a homomorphism of Lie algebras. \medskip 

From claim 2.2.1 it follows that the adjoint representation on $\widetilde{\mathfrak{n}}$ of the root subalgebra 
${\widetilde{\mathfrak{n}}}^{(\widehat{\alpha})}$, where $\widehat{\alpha } \in \Pi $, is completely 
reducible. We will identify the representation space of an irreducible 
summand of this representation with the extremal root chain 
${\mathcal{S}}^{\widehat{\beta }}_{\widehat{\alpha }}$ through 
$\widehat{\beta } \in \widehat{\mathcal{R}}$ in the direction $\widehat{\alpha}$. \medskip 

We begin by recalling the concept of an extremal root chain. 
Let $\widehat{\alpha }, \widehat{\beta } \in \widehat{\mathcal{R}}$ and set $V 
= {\spann }_{\R } \{ \widehat{\alpha } \setrule \, \widehat{\alpha } \in \widehat{\mathcal{R}}  \}$. A collection 
of the form $\{ \widehat{\beta } + j\, \widehat{\alpha } \in V \setrule \,  j \in F \} $, where 
$F$ is a finite subset of $\Z $, is called a \emph{chain}. If a chain is of the form 
\begin{equation}
\widehat{\beta } - q\, \widehat{\alpha} , \, \, \widehat{\beta } - (q-1) \widehat{\alpha} , \ldots , 
\widehat{\beta } - \widehat{\alpha } , \, \, \widehat{\beta }, \, \, \widehat{\beta} + \widehat{\alpha} , \ldots , 
\widehat{\beta} + p\, \widehat{\alpha} ,
\label{eq-sec2subsec3one}
\end{equation}
where $q,p \in {\Z}_{\ge 0}$ and $\ell \in F$ if and only if 
$\ell \in \Z $ and  $-q \le \ell \le p$, then (\ref{eq-sec2subsec3one}) is an 
\emph{unbroken chain through} $\widehat{\beta }$ \emph{in the direction} $\widehat{\alpha }$ with 
\emph{with integer pair} $(q,p)$. Its \emph{length} is $q+p+1$. Given a chain 
$\{ \widehat{\beta} + j\, \widehat{\alpha } \in V\setrule \,  j \in F \} $, 
a chain of the form $\{ \widehat{\beta} + j\, \widehat{\alpha} \in V \setrule \,  j \in F' \subseteq F \} $ is a 
\emph{subchain} of the given chain. This subchain is \emph{proper} if $F'$ is a 
proper subset of $F$. If each element of a chain lies in $\widehat{\mathcal{R}} \cup 
\{ \widehat{0} \} $, where $\widehat{0} = \zeta $, then the chain is called a \emph{root chain} 
through $\widehat{\beta }$ in the direction $\widehat{\alpha }$. If (\ref{eq-sec2subsec3one}) is an unbroken root chain with integer pair $(q,p)$ chosen as large as possible so that this chain is not a proper subchain of another unbroken root chain through $\widehat{\beta }$ in the direction $\widehat{\alpha }$, then the unbroken root chain ${\mathcal{S}}^{\widehat{\beta }}_{\widehat{\alpha }}$ 
(\ref{eq-sec2subsec3one}) is called an \emph{extremal root chain} through $\widehat{\beta }$ in the 
direction $\widehat{\alpha }$ with integer pair $(q,p)$. Since axiom 3 of a system of roots holds, 
for every $\widehat{\beta }$, $\widehat{\alpha} \in 
\widehat{\mathcal{R}}$ there is an extemal roots chain with a suitable integer pair $(q,p)$.  \medskip 

We now identify the irreducible summands of the adjoint representation of 
${\widetilde{\mathfrak{n}}}^{(\widehat{\alpha})}$, $\widehat{\alpha } \in \widehat{\Pi}$, on 
$\widetilde{\mathfrak{n}}$. Corresponding to the extremal root chain 
(\ref{eq-sec2subsec3one}) is the subspace    
\begin{equation}
{\widetilde{\mathfrak{n}}}^{\widehat{\beta }}_{\widehat{\alpha }} = 
\sum _{\stackrel{-q \le j \le p}{\rule{0pt}{6pt}{\scriptscriptstyle j\in \Z }}} \oplus 
{\widetilde{\mathfrak{n}}}_{\widehat{\beta} + j\, \widehat{\alpha }}.  
\label{eq-irred}
\vspace{-.1in} 
\end{equation}
of $\widetilde{\mathfrak{n}}$. If $\widehat{0}$ appears in the root chain (\ref{eq-sec2subsec3one}), then 
${\widetilde{\mathfrak{n}}}_{\widehat{0}} = Z = {\spann }_{\C } \{ X_{\zeta } \} $. If $\widehat{\beta } \pm \widehat{\alpha } \notin \widehat{\mathcal{R}}$, 
then ${\widetilde{\mathfrak{n}}}^{\widehat{\beta }}_{\widehat{\alpha }} = 
{\widetilde{\mathfrak{n}}}_{\widehat{\beta }}$.  
Using the fact that $[{\widetilde{\mathfrak{n}}}_{\widehat{\alpha }}, {\widetilde{\mathfrak{n}}}_{\widehat{\beta }}] = {\widetilde{\mathfrak{n}}}_{\widehat{\alpha } +\widehat{\beta }}$ if $\widehat{\alpha } +\widehat{\beta } \in \widehat{\mathcal{R}}$ and $\{ 0 \} $ if $\widehat{\alpha } +\widehat{\beta } \notin \widehat{\mathcal{R}}$ and the definition of extremal root chain, it follows that 
${\widetilde{\mathfrak{n}}}^{\widehat{\beta }}_{\widehat{\alpha }}$ is an 
${\ad }_{{\widetilde{\mathfrak{n}}}^{(\alpha )}}$-invariant subspace of $\widetilde{\mathfrak{n}}$. \medskip 
 
\noindent \textbf{Claim 2.2.3} For every $\widehat{\alpha} \in \widehat{\Pi}$, $\widehat{\beta} \in \widehat{\mathcal{R}}$ 
the subspace ${\widetilde{\mathfrak{n}}}^{\widehat{\beta }}_{\widehat{\alpha }}$ 
(\ref{eq-irred}) of $\widetilde{\mathfrak{n}}$ is a representation space for the adjoint representation of 
root subalgebra ${\widetilde{\mathfrak{n}}}^{(\widehat{\alpha })}$, which is irreducible. \medskip 

\noindent \textbf{Proof}. Let $(q,p)$ be the integer pair associated 
to the extremal root chain through $\widehat{\beta }$ in the direction of 
$\widehat{\alpha }$. Then ${\widetilde{\mathfrak{n}}}_{\widehat{\beta} + p\, \widehat{\alpha }}$ is the top root space corresponding to the top root $\widehat{\beta} + p\, \widehat{\alpha }$ in this chain. 
The linear operator ${\ad }_{X_{-\widehat{\alpha }}}$ steps down the 
extremal root chain from ${\widetilde{\mathfrak{n}}}_{\widehat{\beta} + j\, \widehat{\alpha }}$ to 
\begin{displaymath}
\left\{ \begin{array}{rl} {\widetilde{\mathfrak{n}}}_{\widehat{\beta} + (j-1)\widehat{\alpha }}, & 
\mbox{when $-q < j \le p$, $j \in \Z$ } \\
\rowspace \{ 0 \} , & \mbox{when $j =-q$.} \end{array} \right. 
\end{displaymath}  
Since ${\ad}_{X_{-\widehat{\alpha }}}{\widetilde{\mathfrak{n}}}_{\widehat{\beta} + q \, \widehat{\alpha }}v = 
\{ 0 \}$, the linear map ${\ad }_{X_{-\widehat{\alpha }}}:
{\widetilde{\mathfrak{n} }}^{\widehat{\beta }}_{\widehat{\alpha }} \rightarrow 
{\widetilde{\mathfrak{n}}}^{\widehat{\beta }}_{\widehat{\alpha }}$ is nilpotent of height $p+q$, because 
${\ad }^{p+q}_{X_{-\widehat{\alpha }}}{\widetilde{\mathfrak{n}}}_{\widehat{\beta } + p\, \widehat{\alpha }} = 
{\widetilde{\mathfrak{n}}}_{\beta + q \, \alpha }$, while ${\ad }^{p+q+1}_{X_{-\widehat{\alpha }}}
{\widetilde{\mathfrak{n}}}_{\widehat{\beta} +  p\, \widehat{\alpha }} = \{ 0 \} $. Taking a nonzero vector 
$v \in {\widetilde{\mathfrak{n}}}_{\widehat{\beta } + p \, \widehat{ \alpha }}$, the above argument shows that 
$\mathfrak{f} = \{ v, {\ad }_{X_{-\widehat{\alpha }}}v, \ldots , 
{\ad }^{p+q}_{X_{-\widehat{\alpha }}}v \} $ is a Jordan chain in 
${\widetilde{\mathfrak{n}}}^{\widehat{\beta }}_{\widehat{\alpha }}$ of length 
$p+q+1$. So $\mathfrak{f}$ is a basis of $W = {\spann}_{\C} \{ v, {\ad }_{X_{-\widehat{\alpha }}}v, \ldots , 
{\ad }^{p+q}_{X_{-{\widehat{\alpha }}}} v \} $ with respect to which 
the matrix of ${\ad }_{X_{-\widehat{\alpha }}}|W$ is a $(p+q+1) \times 
(p+q+1)$ lower Jordan block. By claim 2.2.1 $W$ is a representation 
space for a $p+q+1$-dimensional irreducible representation of 
${\ad }_{{\widetilde{\mathfrak{n}}}^{(\widehat{\alpha })}}$. But 
${\dim }_{\C}\,  {\widetilde{\mathfrak{n}}}_{\widehat{\beta} +j\, \widehat{ \alpha }} = 1$ 
for every $-q \le j \le p$, $j\in \Z$. Therefore 
${\dim }_{\C}\, {\widetilde{\mathfrak{n}}}^{\widehat{\beta }}_{\widehat{\alpha } }= p+q+1 = \dim W$, which shows that $W = {\widetilde{\mathfrak{n}}}^{\widehat{\beta }}_{\widehat{\alpha }}$.  \hfill $\square $ 

\subsection{Verification of axioms 4 and 5}
\label{sec2subsec3}

In this subsection we verify that axioms 4 and 5 for a system of roots holds for the set of roots 
$\widehat{\mathcal{R}}$ of a very special sandwich algebra. \medskip 

First we verify that axiom 5 holds. In other words, \medskip 

\noindent \textbf{Lemma 2.3.1} For every $\widehat{\alpha} \in \widehat{\mathcal{R}}$ the Killing integer 
$\langle \widehat{\alpha }, \widehat{\alpha } \rangle $ associated to the extremal root chain through 
$\widehat{\alpha }$ in the direction of $\widehat{\alpha } $ is $2$. \medskip

\noindent \textbf{Proof}. From lemma 1.2.3 it follows that  
for every $\widehat{\beta }$, $\widehat{\alpha}  \in \widehat{\mathcal{R}}$ we have 
\begin{equation}
\langle \widehat{\beta }, - \widehat{\alpha }\rangle = - \langle \widehat{\beta },\widehat{\alpha } \rangle \, \,  
\mathrm{and}\, \, \langle - \widehat{\beta }, \widehat{\alpha }\rangle = - \langle \widehat{\beta }, \widehat{\alpha} \rangle ,
\label{eq-ktwonew} 
\end{equation}
respectively. Thus it suffices to assume that 
$\widehat{\beta }, \widehat{\alpha }\in \widehat{\Pi }$. We now determine the extremal root 
chain in $\widehat{\mathcal{R}}$ through $\widehat{\alpha} \in \widehat{\Pi }$ in the 
direction $\widehat{\alpha }$. Such an unbroken root chain is 
\begin{equation}
\widehat{\alpha } - 2\, \widehat{\alpha } = -\widehat{\alpha }, \, \, \widehat{\alpha }- 1\widehat{\alpha } = \widehat{0}, 
\, \, \widehat{\alpha } + 0 \, \widehat{\alpha }  =\widehat{\alpha }.   
\label{eq-twodstar}
\end{equation}
The root chain (\ref{eq-twodstar}) is extremal, since $2\widehat{\alpha }  =  
\widehat{\alpha } + \widehat{\alpha } 
\in \mathcal{Z} $, which follows because $X_{2\widehat{\alpha }} \in 
{\widetilde{\mathfrak{n}}}_{2\widehat{\alpha}} = [{\widetilde{\mathfrak{n}}}_{\widehat{\alpha }}, 
{\widetilde{\mathfrak{n}}}_{\widehat{\alpha }}] \in Z$. Thus $2\widehat{\alpha}$ does not lie in 
$\widehat{\mathcal{R}}$. The root chain (\ref{eq-twodstar}) has associated integer pair $(2,0)$. Hence by definition its Killing integer $\langle \alpha , \alpha \rangle$ is $2$. \hfill $\square $ \medskip   

We now verify that axiom 4 of a system of roots holds for the collection of roots $\widehat{\mathcal{R}}$. This amounts to proving  \medskip

\noindent \textbf{Claim 2.3.2} For each fixed $\widehat{\alpha } \in \widehat{\mathcal{R}}$ the function 
\begin{displaymath}
K_{\widehat{\alpha }}: \widehat{\mathcal{R}} \rightarrow \Z: 
\widehat{\beta} \mapsto \langle \widehat{\beta }, \widehat{\alpha} \rangle 
\end{displaymath} 
is linear. By linear we mean: if $\widehat{\gamma } , \widehat{\delta } \in \widehat{\mathcal{R}}$ such that 
$\widehat{\gamma }+ \widehat{\delta } \in \widehat{\mathcal{R}} $, then $K_{\widehat{\alpha }}(\widehat{\gamma } +\widehat{\delta } ) = 
K_{\widehat{\alpha }}(\widehat{\gamma }) + K_{\widehat{\alpha }}(\widehat{\delta })$. \medskip     

We will need a few preliminary results. Let $\widehat{\alpha } \in \widehat{\Pi }$. Suppose that 
${\widehat{\beta }}_1$, ${\widehat{\beta }}_2$ and 
${\widehat{\beta }}_1 + {\widehat{\beta }}_2$ lie in $\widehat{\mathcal{R}}$. For $i=1,2$ let 
\begin{equation}
{\widehat{\beta }}_i - q_i \, \widehat{\alpha }, \ldots , \,  {\widehat{\beta }}_i - 1\alpha , 
\, {\widehat{\beta }}_i, \, {\widehat{\beta }}_i + 1\widehat{\alpha }, \ldots , {\widehat{\beta }}_i + p_i \, \widehat{\alpha }
\label{eq-sec2subsec4two}
\end{equation}
be an extremal root chain through ${\widehat{\beta }}_i$ in the direction $\widehat{\alpha }$ with integer pair 
$(q_i,p_i)$. Let 
\begin{equation}
\mbox{\footnotesize ${\widehat{\beta }}_1 + {\widehat{\beta }}  - r \, \widehat{\alpha } , \ldots , 
{\widehat{\beta }}_1+ {\widehat{\beta }} _2 - (r-1) \widehat{\alpha },  \ldots 
{\widehat{\beta }}_1 +  {\widehat{\beta }}_2 +(s-1) \widehat{\alpha } ,  
{\widehat{\beta }}_1 + {\widehat{\beta }}_2 + s\, \widehat{\alpha }$}
\label{eq-sec2subsec4three}
\end{equation}
be an extremal root chain through ${\widehat{\beta }}_1+ {\widehat{\beta }}_2$ 
in the direction $\widehat{\alpha }$ with integer pair $(r,s)$. Consider the chain 
\begin{equation}
\mbox{\footnotesize ${\widehat{\beta }}_1 + {\widehat{\beta }}_2 - (q_1+q_2) \widehat{\alpha }, \ldots , 
{\widehat{\beta }}_1+ {\widehat{\beta }}_2 -  0\, \widehat{\alpha} , \ldots , 
{\widehat{\beta }_1} + {\widehat{\beta }}_2 + (p_1 +p_2)\widehat{\alpha }$}
\label{eq-sec2subsec4four}
\end{equation}

\noindent \textbf{Lemma 2.3.3} The chain (\ref{eq-sec2subsec4three}) is an extremal root subchain of the unbroken chain 
(\ref{eq-sec2subsec4four}).\medskip

\noindent \textbf{Proof.} The following argument shows that 
the chain (\ref{eq-sec2subsec4four}) is unbroken. Write 
\begin{align}
\mbox{{\footnotesize ${\widehat{\beta }}_2 -  q_2\widehat{\alpha } + {\widehat{\beta }}_1 - q_1 \widehat{\alpha }, \, 
{\widehat{\beta }}_2 - q_2\widehat{\alpha } + {\widehat{\beta }}_1 - (q_1-1) \widehat{\alpha }, \ldots , 
{\widehat{\beta }}_2 - q_2\widehat{\alpha } + {\widehat{\beta }}_1 +  p_1 \widehat{\alpha }$}} & \notag \\
&\hspace{-3in} \mbox{{\footnotesize ${\widehat{\beta }}_2- (q_2-1)\widehat{\alpha } + {\widehat{\beta }}_1 + p_1 \widehat{\alpha }, \ldots , {\widehat{\beta }}_2 + p_2\widehat{\alpha } + {\widehat{\beta }}_1 + p_1 \widehat{\alpha}  $.}}
\label{eq-sec2subsec4five}
\end{align}
Using the fact that the chains in (\ref{eq-sec2subsec4two}) are unbroken, 
it follows that the chain (\ref{eq-sec2subsec4five}) is unbroken. Clearly the 
chains (\ref{eq-sec2subsec4four}) and (\ref{eq-sec2subsec4five}) are equal. \medskip 

To prove that the chain (\ref{eq-sec2subsec4three}) is a subchain of the chain (\ref{eq-sec2subsec4four}) suppose that 
\begin{equation}
q_1+q_2 < r \quad \mathrm{or} \quad p_1+p_2 < s.  
\label{eq-sec2subsec4six}
\end{equation} 
Recall that ${\widetilde{\mathfrak{n}}}^{(\widehat{\alpha })}$ is the root subalgebra of 
$\widetilde{\mathfrak{n}}$ associated to the positive root $\widehat{\alpha }\in \widehat{\Pi }$. 
Then the root subalgebra ${\widetilde{\mathfrak{n}}}^{(\widehat{\alpha })}$ is equal to $\spann \{ X_{\widehat{\alpha } }, X_{-\widehat{\alpha }}, H_{\zeta } \} $, which is isomorphic to ${\mathrm{h}}_3$. Let 
\begin{displaymath}
W = \sum_{\stackrel{-(q_1+q_2) \le \ell \le (p_1+p_2)}{\rule{0pt}{6pt}{\scriptscriptstyle \ell \in \Z}}} \hspace{-15pt} \oplus \, {\widetilde{\mathfrak{n}}}_{{\widehat{\beta }}_1+ {\widehat{\beta }}_2 + \ell \widehat{\alpha }}, 
\end{displaymath}
where ${\widetilde{\mathfrak{n}}}_{\widehat{\gamma }}$ is the root space of $\widetilde{\mathfrak{n}}$ 
corresponding to $\widehat{\gamma }$ in the chain (\ref{eq-sec2subsec4four}). Since the root 
chains in (\ref{eq-sec2subsec4two}) are unbroken and extremal,  
\begin{displaymath}
{\ad }^{q_i+p_i+1}_{X_{-\widehat{\alpha}}} {\widetilde{\mathfrak{n}}}_{{\widehat{\beta }}_i +p_i\, \widehat{\alpha }}= \{ 0 \}  
\end{displaymath}
and therefore 
\begin{displaymath}
{\ad }^{p_1+p_2+q_1+q_2+1}_{X_{-\widehat{\alpha}}} 
{\widetilde{\mathfrak{n}}}_{{\widehat{\beta }}_1 + {\widehat{\beta}}_2 +(p_1+p_2)\widetilde{\alpha }} = \{ 0 \} .  
\end{displaymath}
This implies that  $W$ is ${\ad }_{{\widetilde{\mathfrak{n}}}^{(\widehat{\alpha })}}$-invariant. By claim 2.2.3 
\begin{displaymath}
{\widetilde{\mathfrak{n}}}^{{\widehat{\beta }}_1 +{\widehat{\beta }}_2}_{\widehat{\alpha }} = 
\sum_{\stackrel{-r \le m \le s}{{\scriptscriptstyle m \in \Z}}} \oplus 
\, {\widetilde{\mathfrak{n}}}_{{\widehat{\beta }}_1 + {\widehat{\beta }}_2 + m \widehat{\alpha }}
\vspace{-.1in} 
\end{displaymath} 
is a representation space for an irreducible representation of 
${\ad }_{{\widetilde{\mathfrak{g}}}^{(\widehat{\alpha })}}$ on 
${\widetilde{\mathfrak{n}}}^{{\widehat{\beta }}_1 +{\widehat{\beta }}_2}_{\widehat{\alpha }}$, 
since (\ref{eq-sec2subsec4three}) is an extremal root chain. From the hypothesis (\ref{eq-sec2subsec4six}) it follows that $W$ is a proper ${\ad }_{{\widetilde{\mathfrak{n}}}^{(\widehat{\alpha })}}$-invariant 
subspace of ${\widetilde{\mathfrak{n}}}^{{\widehat{\beta }}_1 + {\widehat{\beta }}_2}_{\widehat{\alpha }}$. 
(Recall that equation (\ref{eq-irred}) implies that ${\widetilde{\mathfrak{n}}}_{\widehat{\gamma }}  = \{ 0 \} $, if 
$\widehat{\gamma }= {\widehat{\beta }}_1+{\widehat{\beta }}_2 + m \, \widehat{\alpha }$ lies in the chain 
(\ref{eq-sec2subsec4four}) and either $m< -r$ or $m >s$.) But this contradicts the irreducibility of the 
${\ad }_{{\widetilde{\mathfrak{n}}}^{(\widehat{\alpha })}}$-representation on 
${\widetilde{\mathfrak{n}}}^{{\widehat{\beta }}_1 + {\widehat{\beta }}_2}_{\widehat{\alpha }}$. Therefore 
the statement in (\ref{eq-sec2subsec4six}) is false, that is, 
\begin{equation}
q_1+q_2 \ge r \quad \mathrm{and} \quad p_1+p_2 \ge s . 
\label{eq-sec2subsec4seven}
\end{equation}
This shows that the chain (\ref{eq-sec2subsec4three}) is a subchain of the chain (\ref{eq-sec2subsec4four}). Since the chain (\ref{eq-sec2subsec4three}) is extremal by hypothesis, we have proved the lemma. \hfill $\square $ \medskip 

We now define the notion of the sum of the extremal root chains 
\begin{equation}
{\widehat{\beta }}_i - q_i\, \widehat{\alpha} , \ldots , {\widehat{\beta }}_i + p_i \, \widehat{\alpha} 
\label{eq-sec2subsec4eight}
\end{equation}
with integer pair $(q_i,p_i)$ for $i=1,2$. Suppose that ${\widehat{\beta }}_1 + {\widehat{\beta }}_2 \in 
\widehat{\mathcal{R}}$. The \emph{sum} of the extremal root chains (\ref{eq-sec2subsec4eight}) is the root chain obtained from the unbroken chain 
\begin{equation}
{\widehat{\beta }}_1 + {\widehat{\beta }}_2  - (q_1+q_2)\widehat{\alpha }, \ldots , 
{\widehat{\beta }}_1 + {\widehat{\beta }}_2 + (p_1+p_2)\widehat{\alpha }
\label{eq-sec2subsec4nine}
\end{equation}
by removing all the elements ${\widehat{\beta }}_1 + {\widehat{\beta }}_2 + j \, \widehat{\alpha }$ of 
(\ref{eq-sec2subsec4nine}) which do not lie in $\widehat{\mathcal{R}} \cup \{ \widehat{0} \}$, 
where $\widehat{0} = \zeta$. By lemma 2.3.3 the extremal root chain  
\begin{equation}
{\widehat{\beta }}_1 + {\widehat{\beta }}_2 - r\, \widehat{\alpha} , \ldots , 
{\widehat{\beta }}_1 + {\widehat{\beta }}_2 + s \, \widehat{\alpha} 
\label{eq-sec2subsec4ten}
\end{equation}
with integer pair $(r,s)$ is an unbroken root subchain of (\ref{eq-sec2subsec4nine}). From the definition of 
the nonnegative integer $r$ it follows that ${\widehat{\beta }}_1 + {\widehat{\beta }}_2 - k \, \widehat{\alpha } 
\notin \widehat{\mathcal{R}} \cup \{ \widehat{0} \} $ for every $k \in \Z$ with $r+1 \le k \le (q_1+q_2)$. 
Similarly, from the definition of the nonnegative integer $s$ it follows that 
$\{ {\widehat{\beta }}_1 + {\widehat{\beta }}_2 +  m\, \widehat{\alpha }
\notin \widehat{\mathcal{R}} \cup \{ \widehat{0} \} \} $ for every $m \in \Z$ with $s+1 \le m \le (p_1+p_2)$. Consequently, the sum of the extremal root chains (\ref{eq-sec2subsec4eight}) is the extremal root chain (\ref{eq-sec2subsec4ten}).  \medskip 

Next we show that \medskip 

\noindent \textbf{Claim 2.3.4} The integer $(q_1+q_2)-(p_1+p_2)$ associated to the unbroken chain 
(\ref{eq-sec2subsec4nine}) with integer pair $(q_1+q_2,p_1+p_2)$ is equal to the Killing integer $r-s$ 
of the extremal root chain (\ref{eq-sec2subsec4ten}).  \medskip 

\noindent \textbf{Proof.} For $i=1,2$ let $\widehat{\alpha} \in \widehat{\Pi }$ and ${\widehat{\beta }}_i \in \widehat{\mathcal{R}}$. 
Recall that  ${\widehat{\beta }}_i - q_i \, \widehat{\alpha } , \ldots , {\widehat{\beta }}_i +  p_i\, \widehat{\alpha }$ is an extremal root chain with integer pair $(q_i,p_i)$ and length $q_i+p_i+1$. Since (\ref{eq-sec2subsec4eight}) is an extremal root chain, 
${\widetilde{\mathfrak{n}}}^{{\widehat{\beta }}_i}_{\widehat{\alpha }}$ is a complex $N_i=q_i+p_i+1$-dimensional representation space for an irreducible ${\ad }_{{\widetilde{\mathfrak{n}}}^{(\widehat{\alpha })}}$ representation. Let $n_i =${\scriptsize $ \left\{ \begin{array}{rl} \onehalf (q_i+p_i), & 
\mbox{if $N_i$ is odd} \\
\onehalf (q_i+p_i+1), & \mbox{if $N_i$ is even.} \end{array} \right. $} 
Then $n_i$ is a nonnegative integer. \medskip 

Recall that the root subalgebra ${\widetilde{\mathfrak{n}}}^{(\widehat{\alpha })}$ is isomorphic to 
${\mathrm{h}}_3$. Let $\iota :\widetilde{\mathfrak{n}} \rightarrow \underline{\mathfrak{g}}$ 
be the inclusion map, which is defined since $\widetilde{\mathfrak{n}}$ is a 
subalgebra of the \linebreak 
simple Lie algebra $\underline{\mathfrak{g}}$. By claim 2.1.2 the image of the root vector $X_{-\widehat{\alpha }}$ under $\iota $ is the root vector $X_{-\underline{\alpha }}$ in 
$\underline{\mathfrak{g}}$. Since $\underline{\mathfrak{g}}$ is simple, this latter root vector 
embeds into the root subalgebra ${\underline{\mathfrak{g}}}^{(\underline{\alpha })} = 
\spann \{ X_{-\underline{\alpha }}, X_{\underline{\alpha }}, H_{\underline{\alpha }} \} $ of 
$\underline{\mathfrak{g}}$, which is isomorphic to $\ssl (2, \C)$. Since the root chains in 
(\ref{eq-sec2subsec4eight}) are extremal, it follows 
that the complex $N_i =q_i+p_i+1$-dimensional ${\ad }_{{\widetilde{\mathfrak{n}}}^{(\widehat{\alpha })}}$ representation on ${\widetilde{\mathfrak{n}}}^{{\widehat{\beta }}_i}_{\widehat{\alpha}}$ is irreducible. Thus there is a vector $v_i \in {\widetilde{\mathfrak{n}}}^{{\widehat{\beta }}_i}_{\widehat{\alpha}}$ such that $V_i = 
\spann \{ v_i, {\ad }_{X_{-\widehat{\alpha }}}v_i, \ldots , {\ad }^{q_i+p_i}_{X_{-{\widehat{\alpha }}}} v_i \} = 
{\widetilde{\mathfrak{n}}}^{{\widehat{\beta }}_i}_{\widehat{\alpha}}$. Consider the vector $w_i = \iota (v_i) \in \underline{\mathfrak{g}}$. Let $W_i = \spann \{ w_i, {\ad }_{X_{-\underline{\alpha }}}w_i, \ldots ,$ ${\ad }^{q_i+p_i}_{X_{-\underline{\alpha }}}w_i \} $. Because  
$X_{-\underline{\alpha }} = \iota (X_{-\widehat{\alpha }})$ and ${\ad }_{X_{-\underline{\alpha }}} \comp \iota = 
\iota \comp {\ad }_{X_{-\widehat{\alpha}}}$, it follows that $W_i = \iota (V_i)$ and 
${\ad }^{q_i+p_i+1}_{X_{-\underline{\alpha }}}w_i =0$. Since $V_i$ is an $N_i$-dimensional irreducible representation of ${\widetilde{\mathfrak{n}}}^{(\widehat{\alpha })}$ 
on $\widetilde{\mathfrak{n}}$, we find that the root subalgebra 
${\underline{\mathfrak{g}}}^{(\underline{\alpha})} $ of $\underline{\mathfrak{g}}$, 
which is isomorphic to $\ssl (2, \C)$ and contains $\iota ({\widetilde{\mathfrak{n}}}^{(\widehat{\alpha } )})$, 
acts irreducibly on $W_i$. Therefore ${\ad }_{H_{\underline{\alpha }}}w_i = n_i w_i$, 
where $n_i$ is a nonnegative integer. This shows that 
\begin{equation}
-n_i, \, \, -n_i+2, \, \, -n_i+2(2), \ldots , -n_i+2(n_i-1) = 
n_i -2, \, \, n_i 
\label{eq-sec2subsec4twelveb}
\end{equation}
lists the elements of the extremal root chains in (\ref{eq-sec2subsec4eight}). The root chain (\ref{eq-sec2subsec4ten})  
with integer pair $(r,s)$ and length $M =r+s+1$ is extremal. Therefore 
${\widetilde{\mathfrak{n}}}^{{\widehat{\beta }}_1 + {\widehat{\beta }}_2}_{\widehat{\alpha }}$ is an $M$-dimensional representation space for an ${\ad }_{{\widetilde{\mathfrak{n}}}^{(\widehat{\alpha })}}$-irreducible representation. The eigenvalues of ${\ad }_{H_{\underline{\alpha }}}$ on 
$i({\widetilde{\mathfrak{n}}}^{{\widehat{\beta }}_1 +{\widehat{\beta }}_2}_{\widehat{\alpha }})$ are listed as follows.
\begin{equation}
-m, \, \, -m+2, \, \, -m+2(2), \ldots , -m+2(m-1) = m-2, \, \, m . 
\label{eq-sec2subsec4elevenstarf}
\end{equation}
Since the vector spaces 
$i({\widetilde{\mathfrak{n}}}^{{\widehat{\beta }}_1+ {\widehat{\beta }}_2}_{\widehat{\alpha }})$ 
and ${\underline{\mathfrak{g}}}^{{\underline{\beta }}_1 + {\underline{\beta }}_2}_{\underline{\alpha }}$
are equal, equation (\ref{eq-sec2subsec4elevenstarf}) lists the elements of the extremal root chain 
(\ref{eq-sec2subsec4ten}) provided that the positive integer $m=${\scriptsize $\left\{ \begin{array}{rl} 
\onehalf (r+s), & \mbox{if $M$ is odd} \\ \onehalf (r+s+1), & 
\mbox{if $M$ is even} \end{array} \right. $.} \medskip 

From claim 2.3.3 the extremal root chain (\ref{eq-sec2subsec4ten}) is a subchain 
of the unbroken chain (\ref{eq-sec2subsec4nine}). Thus the list (\ref{eq-sec2subsec4elevenstarf}) 
is a sublist of the list 
\begin{align}
&-\ell =-(n_1+n_2), \, \, -\ell +2, \, \, -\ell +2(2), \ldots , \notag \\
&\hspace{1in} -\ell +2(n_1+n_2-1) = \ell -2, \, \, \ell ,  
\label{eq-sec2subsec4thirteen}
\end{align}
which labels the elements of the unbroken chain (\ref{eq-sec2subsec4nine}). Suppose 
that the nonnegative integers $\ell $ and $m$ have a different parity. In particular, suppose that $\ell $ is even and $m$ is odd. Then $1$ appears in the list (\ref{eq-sec2subsec4elevenstarf}) but not in the list (\ref{eq-sec2subsec4thirteen}). This contradicts the fact that 
(\ref{eq-sec2subsec4elevenstarf}) is a sublist of (\ref{eq-sec2subsec4thirteen}). Suppose that $\ell $ is odd 
and $m$ is even. Then $0$ appears in the list (\ref{eq-sec2subsec4elevenstarf}) but not in the list (\ref{eq-sec2subsec4thirteen}). Again this contradicts the fact that (\ref{eq-sec2subsec4elevenstarf}) is a sublist of (\ref{eq-sec2subsec4thirteen}). Therefore 
$\ell $ and $m$ must have the same parity. Conversely, if $\ell $ and $m$ have the 
same parity, then the list (\ref{eq-sec2subsec4elevenstarf}) is a sublist of the list 
(\ref{eq-sec2subsec4thirteen}). Consequently, for $i=1,2$ there is a $j_i \in \Z$, $0 \le j_i \le \ell $, such that 
\begin{displaymath}
-\ell +2j_1 = -m \quad \mathrm{and} \quad \ell -2j_2 =m.
\end{displaymath}
Therefore $j_1=j_2 =j$. Thus $j$ is the number of elements of the 
unbroken chain (\ref{eq-sec2subsec4nine}) which need to be removed from its left and right ends to obtain the extremal root chain (\ref{eq-sec2subsec4ten}). Hence the integer $(q_1+q_2) - (p_1+p_2)$ associated to the unbroken chain (\ref{eq-sec2subsec4nine}) with integer pair $(q_1+q_2, p_1+p_2)$ is equal to the Killing 
integer $r-s$ of the extremal root chain (\ref{eq-sec2subsec4ten}) with integer pair $(r,s)$. This proves claim 2.3.4. \hfill $\square $ \medskip  

Claim 2.3.4 may be reformulated as \medskip 

\noindent \textbf{Claim 2.3.5} Let $\widehat{\alpha } \in \widehat{\Pi }$ with ${\widehat{\beta }}_1$, and ${\widehat{\beta }}_2 \in \widehat{\mathcal{R}}$ such that ${\widehat{\beta }}_1 + {\widehat{\beta }} _2 \in \widehat{\mathcal{R}}$. For $i=1,2$ let 
\begin{equation}
{\widehat{\beta }}_i -  q_i \, \widehat{\alpha }, \ldots , {\widehat{\beta }}_i  +  p_i \, \widehat{\alpha }
\label{eq-sec2subsec4twelvec}
\end{equation}
be extremal root chains with integer pair $(q_i,p_i)$ and Killing integer 
$\langle {\widehat{\beta }}_i, \widehat{\alpha } \rangle = q_i-p_i$. Suppose that the extremal root chain 
\begin{displaymath}
{\widehat{\beta }}_1+ {\widehat{\beta}}_2 - r \, \widehat{\alpha}, \ldots , {\widehat{\beta }}_1+ {\widehat{\beta }}_2 + s \, \widehat{\alpha}
\end{displaymath}
with integer pair $(r,s)$ and Killing integer $\langle {\widehat{\beta }}_1 + {\widehat{\beta }}_2, \widehat{\alpha } \rangle = r-s$ is the sum of the extremal root chains in (\ref{eq-sec2subsec4twelvec}). Then 
\begin{equation}
\langle {\widehat{\beta }}_1+ {\widehat{\beta }}_2, \widehat{\alpha} \rangle = 
\langle {\widehat{\beta }}_1 , \widehat{\alpha } \rangle + 
\langle {\widehat{\beta }}_2 , \widehat{\alpha } \rangle . 
\label{eq-sec2subsec4thirteend}
\end{equation}

\noindent \textbf{Proof.} 
If $-\widehat{\beta } \in \widehat{\Pi }$ and $\widehat{\alpha } \in \widehat{\Pi }$, then by lemma 1.2.3 we have 
$\langle \widehat{\beta }, \widehat{\alpha } \rangle = - \langle -\widehat{\beta }, \widehat{\alpha} \rangle 
=-K_{\widehat{\alpha }}(-\widehat{\beta })$. Similarly, if $\widehat{\beta }  \in \widehat{\Pi }$ and $\widehat{\alpha } \in 
\widehat{\Pi }$, then by lemma 1.2.3 $\langle \widehat{\beta }, \widehat{\alpha }\rangle 
= - \langle \widehat{\beta }, - \widehat{\alpha} \rangle = -K_{-\widehat{\alpha }}(\widehat{\beta })$. Therefore it suffices to prove the claim 
when $\widehat{\alpha}$, $\widehat{\beta } \in \widehat{\Pi }$. 
Using claim 2.3.4 we deduce that for every root $\widehat{\alpha }$ in $\widehat{\Pi }$ the function 
$K_{\widehat{\alpha }}:\widehat{\mathcal{R}} \rightarrow \Z: \widehat{\beta  }\mapsto 
\langle \widehat{\beta }, \widehat{\alpha }\rangle $ is linear.  \hfill $\square $ \medskip 

This completes the verification of \medskip 

\noindent \textbf{Theorem 2.3.6} Let $\widetilde{\mathfrak{g}} = \mathfrak{g} \oplus \widetilde{\mathfrak{n}}$ be a fundamental sandwich algebra. Then the 
collection of roots $\widehat{\mathcal{R}}$ on a Cartan subalgebra $\mathfrak{h}$ of
$\mathfrak{g}$ acting on the nilradical 
$\widetilde{\mathfrak{n}}$ of $\widetilde{\mathfrak{g}}$ is a system of roots. \medskip 

Because fundamental sandwich algebras are Lie algebras which are not semisimple and have 
a system of roots, the notion of a system of roots is a conservative generalization of the concept of a root system of a semisimple Lie algeba. \bigskip 

\par \noindent {\Large \bf Appendix}
\label{appendix} \bigskip 

\noindent  In this appendix we show that each of the algebras ${\widetilde{\mathbf{C}}}_{\ell +1}$, 
${\widetilde{\mathbf{G}}}_2$, and ${\widetilde{\mathbf{E}}}_7$ is a fundamental sandwich algebra. \medskip 

The following list gives $\widetilde{\mathfrak{g}}$, $\mathfrak{g}$, $\underline{\mathfrak{g}}$, the Cartan subalgebra $\mathfrak{h}$ of $\mathfrak{g}$, and the roots $\widehat{\mathcal{R}}$ of the adjoint action of 
$\mathfrak{h}$ on $\widetilde{\mathfrak{n}}$. By ${\underline{\mathbf{X}}}^{(\ell )}_n$ we mean the simple Lie algebra of Cartan type $\underline{\mathbf{X}}$ of rank $n-1$ whose Dynkin diagram is obtained by removing the node numbered $\ell $ from the Dykin diagram of the simple Lie algebra of Cartan type $\mathbf{X}$ of rank $n$. We use the following notation. Let ${\{ {\varepsilon }_i \} }^n_{i=1}$ be the 
standard basis for $({\C }^n)^{\ast }$ and let ${\{ e_i \} }^n_{i=1}$ be the standard dual basis of ${\C }^n$. 
${\mathrm{h}}_{2n+1}$ denotes the the Lie algebra of the $2n+1$ dimensional Heisenberg group. \bigskip
\begin{center}
\begin{tabular}{ll}
\rule{.25in}{0in} $1$. & ${\widetilde{\mathbf{C}}}_{\ell +1}$, \, ${\mathbf{C}}_{\ell } = 
{\mathbf{C}}^{(1)}_{\ell +1}$, \, 
${\underline{\mathbf{C}}}_{\ell +1}$; \, 
$\widetilde{\mathfrak{n}} = {\mathrm{h}}_{2\ell +1}$; \\
\rule{0pt}{12pt} & $\mathfrak{h}$ with basis $\big\{ e_i-e_{i+1}, \, 1 \le i \le \ell \big\} $, \, 
$\underline{\widetilde{H}} = e_{\ell +1}$; \\
\rule{0pt}{26pt}& $\zeta  = -2{\varepsilon }_{\ell +1}|\mathfrak{h}$, \quad  $\widehat{\mathcal{R}} = 
\left\{ \begin{array}{rl}
{\widehat{\alpha }}_k & =  ({\varepsilon }_k-{\varepsilon }_{\ell +1})|\mathfrak{h}, \\ 
{\widehat{\alpha }}_{\ell + k} & = -{\varepsilon }_k-{\varepsilon }_{\ell +1})|\mathfrak{h}, \, 1 \le k \le \ell .
\end{array} \right. $ \vspace{.075in} 
\end{tabular}  
\medskip 

\noindent \begin{tabular}{ll}
$2. $ &${\widetilde{\mathbf{G}}}_2$, \, ${\mathbf{A}}_1 = {\mathbf{G}}^{(2)}_2$, \, ${\underline{\mathbf{G}}}_2$; \, $\widetilde{\mathfrak{n}} = {\mathrm{h}}_5 $; \\ 
& $\mathfrak{h}$ with basis $\big\{ e_1-e_2 \big\} $, \, \, 
$\underline{\widetilde{H}} = -e_1-e_2+2e_3$ \\
\rule{0pt}{15pt}& $\zeta  = ({\varepsilon }_1+{\varepsilon }_2-2{\varepsilon }_3)|\mathfrak{h}$, \\
\rule{0pt}{20pt} & $\widehat{\mathcal{R}} = \left\{ \begin{array}{rl} 
{\widehat{\alpha }}_1 & = (-{\varepsilon }_1-{\varepsilon}_2)|\mathfrak{h}, \, \, 
{\widehat{\alpha }}_2  = ({\varepsilon }_2 -{\varepsilon }_3)|\mathfrak{h}  \\
{\widehat{\alpha }}_3 & = (2{\varepsilon }_1- {\varepsilon }_2-{\varepsilon }_3)|\mathfrak{h}, \, \,  
{\widehat{\alpha }}_4   = 
({\varepsilon }_1 +2{\varepsilon }_2-{\varepsilon }_3)|\mathfrak{h} \end{array} \right. $ 
\end{tabular}  \medskip
 
\noindent \begin{tabular}{ll}
\rule{.35in}{0in}$3.$ &${\widetilde{\mathbf{E}}}_7$, \, ${\mathbf{D}}_6 = {\mathbf{E}}^{(1)}_7$, 
\, ${\underline{\mathbf{E}}}_7$; \, 
$\widetilde{\mathfrak{n}} = {\mathrm{h}}_{33}$; \\ 
& $\mathfrak{h}$ with basis $\big\{ e_2 \pm e_1; \,  
e_{2+i}-e_{1+i}, \, 1 \le i \le 4; \, e_7 \big\} $, \, \, $\underline{\widetilde{H}} = e_7$  \\
\rule{0pt}{15pt} & $\zeta  = -{\varepsilon }_7|\mathfrak{h}$, \\ 

\rule{0pt}{70pt} & $\widehat{\mathcal{R}} = \left\{ \begin{array}{rl} 
{\widehat{\alpha }}_{10}  & = 
\onehalf ( \sum^6_{i=1}{\varepsilon}_i - {\varepsilon }_7)|\mathfrak{h} \\  
\rule{0pt}{12pt} {\widehat{\alpha }}^{\, 10} & = -\onehalf (\sum^7_{i=1}{\varepsilon }_i)|\mathfrak{h},  \\
&\rule{0pt}{12pt} \mbox{For $j$, $\ell \in \{ 1, \ldots ,6\}$ with $j < \ell $}  \\
{\widehat{\alpha }}_{j\ell } & = \onehalf ( \sum^6_{i=1}(-1)^{k(i)}{\varepsilon}_i - {\varepsilon }_7)|\mathfrak{h} \\ 
&\hspace{.25in} \parbox[t]{3.5in}{with $k(j) = k(\ell ) =1$; $k(i) =0$ \\ for $i \in \{ 1, \ldots , 6\} $ where 
$i \ne j$ and $i \ne \ell $}  \\
\rule{0pt}{12pt} {\widehat{\alpha }}^{j\ell } & = \onehalf ( \sum^6_{i=1}(-1)^{k(i)}{\varepsilon}_i - {\varepsilon }_7)|\mathfrak{h} \\  
&\hspace{.25in} \parbox[t]{3.5in}{with $k(j) = k(\ell ) =0$; $k(i) =1$ \\ for $i \in \{1, \ldots , 6\}$ where 
$i \ne j$ and $i \ne \ell $.} \end{array} \right. $ 
\end{tabular} 
\end{center} 

\noindent Table 1. System of roots for fundamental very special sandwich algebras. 
\clearpage 

We verify the entries in table 1 case by case. \medskip 

\noindent 1. ${\widetilde{\mathbf{C}}}_{\ell}$.
\begin{displaymath}
\underline{\mathcal{R}} = \{ \pm ({\varepsilon }_i \pm {\varepsilon }_j), \, 1 \le i < j \le \ell +1; \, 
\pm 2 {\varepsilon }_i, \, 1 \le i \le \ell +1 \} 
\end{displaymath}
is a root system for $\underline{\mathfrak{g}}= {\mathbf{C}}_{\ell +1}$ corresponding to the Cartan 
subalgebra $\underline{\mathfrak{h}}$ with basis $\{ e_i - e_{i+1}, \, 1 \le i \le \ell ; \, 2e_{\ell +1} \} $. 
Let $\underline{\widetilde{H}} = 2e_{\ell +1} \in \underline{\mathfrak{h}}$. Then 
\begin{displaymath}
\mathcal{R} = \{ \underline{\alpha } \in \underline{\mathcal{R}} \setrule \, 
\underline{\alpha }(\underline{\widetilde{H}}) =0 \} = 
\{ \pm ({\varepsilon }_i \pm {\varepsilon }_j), \, 1 \le i < j \le \ell ;\, \pm 2{\varepsilon }_i, \, 1 \le i \le \ell \} 
\end{displaymath}
is the root system of $\mathfrak{g} = {\mathbf{C}}_{\ell } = {\mathbf{C}}^{(1)}_{\ell +1}$ corresponding 
to the Cartan subalgebra $\mathfrak{h}$ with basis  $\{ e_i -e_{i+1}, \, 1 \le i \le \ell -1; \, 2e_{\ell } \} $. Let 
\begin{displaymath}
{\mathcal{R}}^{-} = \{ \underline{\alpha } \in \underline{\mathcal{R}} \setrule \, 
\underline{\alpha }(\underline{\widetilde{H}}) < 0 \} = 
\{ -2{\varepsilon }_{\ell +1}; \, \pm {\varepsilon }_i - {\varepsilon }_{\ell +1}, \, 1\le i \le \ell \} . 
\end{displaymath}
Then 
\begin{displaymath} 
\mbox{\begin{tabular}{l}
$\widehat{0} = \zeta  = -2{\varepsilon }_{\ell +1}|\mathfrak{h}$, \\
\rule{0pt}{20pt}$\widehat{\mathcal{R}} = 
\left\{ \begin{array}{rl}
{\widehat{\alpha }}_k & =  ({\varepsilon }_k-{\varepsilon }_{\ell +1})|\mathfrak{h}, \\ 
{\widehat{\alpha }}_{\ell + k} & = (-{\varepsilon }_k-{\varepsilon }_{\ell +1})|\mathfrak{h}, \, 1 \le k \le \ell .
\end{array} \right. $  
\end{tabular} }
\end{displaymath} 
is the set $\mathfrak{R}$ of roots whose root vectors $X_{\widehat{0}}$, $X_{\widehat{\alpha }}$, 
$\widehat{\alpha} \in \widehat{\mathcal{R}}$ in $\underline{\mathfrak{g}}$ span a complex vector 
space $\widetilde{\mathfrak{n}}$. \medskip 

We now prove some properties of $\widehat{\mathcal{R}}$, which we use in claim A2 to show that 
$\widetilde{\mathfrak{n}}$ is a sandwich algebra whose center is spanned by $X_{\widehat{0}}$. 
For every $1 \le p \le 2\ell $
\begin{equation}
{\widehat{\alpha }}_p + \widehat{0} \notin \widehat{\mathcal{R}}, 
\label{eq-Aone}
\end{equation}
since 
\begin{displaymath}
-2{\varepsilon }_{\ell +1} + \left\{ \begin{array}{rl} {\varepsilon }_p -{\varepsilon }_{\ell +1}, & 
\mbox{if $1 \le p \le \ell $} \\
-{\varepsilon }_{p-\ell } - {\varepsilon }_{\ell +1}, & \mbox{if $\ell +1 \le p \le 2\ell $} \end{array} \right. 
\notin \underline{\mathcal{R}}. 
\end{displaymath}
For every $1 \le j \le \ell $
\begin{equation}
{\widehat{\alpha }}_j + {\widehat{\alpha }}_{\ell +j} = \widehat{0},  
\label{eq-Atwo}
\end{equation}
since 
\begin{displaymath}
({\varepsilon}_j - {\varepsilon}_{\ell +1}) + (-{\varepsilon}_j - {\varepsilon}_{\ell +1}) = -2{\varepsilon }_{\ell +1} 
\in {\mathcal{R}}^{-}. 
\end{displaymath}
For evey $1 \le j \le \ell $ and $\ell +1 \le k \le 2\ell $ with $k \ne j + \ell $ 
\begin{equation}
{\widehat{\alpha }}_j + {\widehat{\alpha }}_k \notin \widehat{\mathcal{R}}, 
\label{eq-Athree}
\end{equation}
because 
\begin{displaymath}
({\varepsilon }_j - {\varepsilon }_{\ell +1}) + (-{\varepsilon }_{k-\ell } - {\varepsilon }_{\ell +1}) \notin 
\widehat{\mathcal{R}}. 
\end{displaymath}
Thus the set of roots $\widehat{\mathcal{R}}$ is closed under addition. Using claim A1 it follows 
that $\widetilde{\mathfrak{n}}$ is the nilradical of the Lie subalgebra $\widetilde{\mathfrak{g}} = 
\mathfrak{g} \oplus \widetilde{\mathfrak{n}}$. Moreover, $\widetilde{\mathfrak{n}}$ is an ideal 
of $\underline{\mathfrak{g}}$. 

\noindent 2. ${\widetilde{\mathbf{G}}}_2$.
\begin{align*}
\underline{\mathcal{R}} & = \pm \left\{ {\varepsilon}_1 - {\varepsilon}_2, \, {\varepsilon}_1 - {\varepsilon}_3, \, 
{\varepsilon}_2 - {\varepsilon}_3 \right.  \\
&\hspace{.25in} \left. 2{\varepsilon}_1-{\varepsilon }_2-{\varepsilon}_3, \, 
-{\varepsilon}_1 + 2{\varepsilon}_2-{\varepsilon }_3, \, 
-{\varepsilon}_1 -{\varepsilon}_2 +2{\varepsilon}_3 \right\} 
\end{align*}
is the root system for $\underline{\mathfrak{g}} = {\mathbf{G}}_2$ corresponding to the Cartan 
subalgebra $\underline{\mathfrak{h}}$ with basis $\{ e_1-e_2, \, e_1 - e_3, \, e_2-e_3 \}$. 
Let $\underline{\widetilde{H}} = -e_1-e_2+2e_3$. Then 
\begin{displaymath}
\mathcal{R} = \{ \underline{\alpha } \in \underline{\mathcal{R}} \setrule \, 
\underline{\alpha }( -e_1-e_2+2e_3) = 0 \} = \{ \pm ({\varepsilon}_1-{\varepsilon}_2) \} 
\end{displaymath}
is the root system of $\mathfrak{g} = {\mathbf{A}}_1 = {\mathbf{G}}^{(2)}_2$ corresponding to the 
Cartan subalgebra $\mathfrak{h}$ with basis $\{ e_1-e_2 \}$. Let 
\begin{align*}
{\mathcal{R}}^{-} & = \{ \underline{\alpha } \in \underline{\mathcal{R}} \setrule \, 
\underline{\alpha}( -e_1-e_2+2e_3) < 0 \} \\
& = \left\{ {\varepsilon}_1 +{\varepsilon}_2 -2{\varepsilon}_3, \, {\varepsilon}_1 - {\varepsilon}_3, \, 
{\varepsilon}_2 - {\varepsilon}_3, \, 2{\varepsilon}_1-{\varepsilon }_2-{\varepsilon}_3, \, 
-{\varepsilon}_1 + 2{\varepsilon}_2-{\varepsilon }_3 \right\} .
\end{align*}
Then $\widehat{0} = \zeta = ({\varepsilon}_1 +{\varepsilon}_2 -2{\varepsilon}_3)|\mathfrak{h}$ and 
\begin{displaymath}
\widehat{\mathcal{R}} = \left\{ \begin{array}{rl} 
{\widehat{\alpha }}_1 = ({\varepsilon}_1 - {\varepsilon}_3)|\mathfrak{h}, & 
{\widehat{\alpha }}_2 = ({\varepsilon}_2 - {\varepsilon}_3)|\mathfrak{h} \\
\rule{0pt}{12pt}{\widehat{\alpha }}_3 = (2{\varepsilon}_1-{\varepsilon }_2-{\varepsilon}_3)|\mathfrak{h}, & 
{\widehat{\alpha }}_4 = (-{\varepsilon}_1 + 2{\varepsilon}_2-{\varepsilon }_3)|\mathfrak{h} 
\end{array} \right. 
\end{displaymath}
is the set of roots $\mathfrak{R}$, whose root vectors $X_{\widehat{0}}; \, X_{\widehat{\alpha }}, 
\widehat{\alpha } \in \widehat{\mathcal{R}}$ in $\underline{\mathfrak{g}}$ span 
a complex vector space $\widetilde{\mathfrak{n}}$. \medskip  

We now prove some properties of $\widehat{\mathcal{R}}$, which we use in claim A2 to show that 
$\widetilde{\mathfrak{n}}$ is a sandwich algebra whose center is spanned by $X_{\widehat{0}}$. First, 
\begin{equation}
\widehat{0} + {\widehat{\alpha}}_j \notin \widehat{\mathcal{R}} \quad \mbox{for every 
$1 \le j \le 4$}, 
\label{eq-Aoneb}
\end{equation}
because 
\begin{displaymath}
\left\{ \begin{array}{l}
({\varepsilon }_1+{\varepsilon }_2-2{\varepsilon }_3) +({\varepsilon }_1 - 
{\varepsilon }_3) ,   \\
({\varepsilon }_1+{\varepsilon }_2-2{\varepsilon }_3) +({\varepsilon }_2 - 
{\varepsilon }_3), \\
({\varepsilon }_1+{\varepsilon }_2-2{\varepsilon }_3) +(2{\varepsilon }_1 - 
{\varepsilon }_2-{\varepsilon }_3),  \\
({\varepsilon }_1+{\varepsilon }_2-2{\varepsilon }_3) +(-{\varepsilon }_1 + 
2{\varepsilon }_2-{\varepsilon }_3)   
\end{array} \right. 
\end{displaymath} 
do not lie in $\underline{\mathcal{R}}$, respectively. Second,  
\begin{equation}
\left\{ \begin{array}{l}
{\widehat{\alpha}}_1 +  {\widehat{\alpha }}_2 = \widehat{0} \\
{\widehat{\alpha }}_3 +  {\widehat{\alpha }}_4 = \widehat{0}. 
\end{array} \right. 
\label{eq-Atwob}
\end{equation}
The first equality above follows since 
$({\varepsilon }_1 - {\varepsilon }_3) + ({\varepsilon }_2 - 
{\varepsilon }_3) = {\varepsilon }_1 +{\varepsilon }_2 - 2{\varepsilon }_3 \in {\mathcal{R}}^{-}$; 
whereas the second equality follows because $(2{\varepsilon }_1 -{\varepsilon }_2 
- {\varepsilon }_3) + (-{\varepsilon }_1 +2{\varepsilon }_2 - {\varepsilon }_3) = 
{\varepsilon }_1 +{\varepsilon }_2 - 2{\varepsilon }_3 \in {\mathcal{R}}^{-}$. 
Third, 
\begin{equation}
{\widehat{\alpha }}_j +  {\widehat{\alpha }}_k \notin 
\widehat{\mathcal{R}}, \quad \mbox{for every } (j,k) \in J = 
\left\{ \begin{array}{l} 
(1,1), \, (1,3), \, (1,4) \\
(2,2), \, (2,3), \, (2,4) \\
(3,1), \, (3,2), \, (3,4) \\
(4,1), \, (4,2), \, (4,4) \end{array} \right\} ,  
\label{eq-Athreeb}
\end{equation}
because 
\begin{displaymath}
\left\{ \begin{array}{l}
({\varepsilon }_1 - {\varepsilon }_3) + ({\varepsilon }_1 - 
{\varepsilon }_3), \, \,    
({\varepsilon }_1 - {\varepsilon }_3) + (2{\varepsilon }_1 -  
{\varepsilon }_2-{\varepsilon }_3),  \\ 
({\varepsilon }_1 - {\varepsilon }_3) + (-{\varepsilon }_1 + 
2{\varepsilon }_2-{\varepsilon }_3), \, \, 
({\varepsilon }_2 - {\varepsilon }_3) + ({\varepsilon }_2 - 
{\varepsilon }_3),  \\
({\varepsilon }_2 - {\varepsilon }_3) + (2{\varepsilon }_1-{\varepsilon }_2 - 
{\varepsilon }_3)  , \, \, 
({\varepsilon }_2 - {\varepsilon }_3) + (-{\varepsilon }_1+2{\varepsilon }_2 - 
{\varepsilon }_3), \\
(2{\varepsilon }_1-{\varepsilon }_2 - {\varepsilon }_3) + 
(2{\varepsilon }_1-{\varepsilon }_2 - {\varepsilon }_3) , \\  
(-{\varepsilon }_1+2{\varepsilon }_2 - {\varepsilon }_3) + 
(-{\varepsilon }_1+2{\varepsilon }_2 - {\varepsilon }_3)  
\end{array} \right. 
\end{displaymath}
do not lie in $\underline{\mathcal{R}}$, respectively. The other cases follow 
from the commutativity of addition, when defined. Thus the set of roots $\widehat{\mathcal{R}}$ is closed under addition. Using claim A1 it follows that $\widetilde{\mathfrak{n}}$ is the nilradical of the Lie subalgebra $\widetilde{\mathfrak{g}} = \mathfrak{g} \oplus \widetilde{\mathfrak{n}}$. Moreover, $\widetilde{\mathfrak{n}}$ is an ideal of $\underline{\mathfrak{g}}$. 
\bigskip

\noindent ${\widetilde{\mathbf{E}}}_7$.
\begin{align*}
\underline{\mathcal{R}} & = \left\{ 
\pm ({\varepsilon}_i \pm {\varepsilon}_j), \, 1 \le i < j \le 6; \, \pm e_7; \, 
\onehalf ( \sum^6_{i=1}(-1)^{k(i)}{\varepsilon }_i \pm {\varepsilon}_7), \, \right. \\
& \hspace{.5in} \left. \mbox{where $k(i) =0$ or $1$ and $\sum^6_{i=1}k(i) \in 2\Z$} \right\} 
\end{align*}
is a root system of $\mathfrak{g} = {\mathbf{E}}_7$ corresponding to the Cartan subalgebra 
$\underline{\mathfrak{h}}$ with basis $\{ e_2 \pm e_1; \, e_{2+i}-e_{1+i}, \, 1 \le i \le 4; \, e_7 \}$. 
Let $\underline{\widetilde{H}} = e_7$. Then 
\begin{displaymath}
\mathcal{R} = \{ \underline{\alpha } \in \underline{\mathcal{R}} \setrule \, \underline{\alpha }(e_7) = 0 \} 
= \{ \pm ({\varepsilon }_i \pm {\varepsilon}_j), \, 1 \le i < j \le 6 \}
\end{displaymath}
is the root system of $\mathfrak{g} = {\mathbf{D}}_6 = {\mathbf{E}}^{(1)}_7$ corresponding to the 
Cartan subalgebra $\mathfrak{h}$ with basis $\{ e_2\pm e_1; \, e_{2+i}-e_{1+i}, \, 1 \le i \le 4 \}$. Let 
\begin{align*}
{\mathcal{R}}^{-} & = \{ \underline{\alpha } \in \underline{\mathcal{R}} \setrule \, 
\underline{\alpha }(e_7) < 0 \} \\
& = \{ -{\varepsilon}_7; \, \onehalf (\sum^6_{i=1}(-1)^{k(i)}{\varepsilon}_i - {\varepsilon}_7), \, 
\mbox{where $k(i) =0$ or $1$ and $\sum^6_{i=1}k(i) \in 2\Z$} \} . 
\end{align*}
Then $\widehat{0} = \zeta = -{\varepsilon}_7|\mathfrak{h}$ and 
\begin{displaymath}
\widehat{\mathcal{R}} = \left\{ \begin{array}{rl} 
{\widehat{\alpha }}_{10}  & = 
\onehalf ( \sum^6_{i=1}{\varepsilon}_i - {\varepsilon }_7)|\mathfrak{h} \\  
\rule{0pt}{12pt} {\widehat{\alpha }}^{\, 10} & = -\onehalf (\sum^7_{i=1}{\varepsilon }_i)|\mathfrak{h},  \\
&\rule{0pt}{12pt} \mbox{For $j$, $\ell \in \{ 1, \ldots ,6\}$ with $j < \ell $}  \\
{\widehat{\alpha }}_{j\ell } & = \onehalf ( \sum^6_{i=1}(-1)^{k(i)}{\varepsilon}_i - {\varepsilon }_7)|\mathfrak{h} \\ 
&\hspace{.25in} \parbox[t]{3.5in}{with $k(j) = k(\ell ) =1$; $k(i) =0$ \\ for $i \in \{ 1, \ldots , 6\} $ where 
$i \ne j$ and $i \ne \ell $}  \\
\rule{0pt}{12pt} {\widehat{\alpha }}^{j\ell } & = \onehalf ( \sum^6_{i=1}(-1)^{k(i)}{\varepsilon}_i - {\varepsilon }_7)|\mathfrak{h} \\  
&\hspace{.25in} \parbox[t]{3.5in}{with $k(j) = k(\ell ) =0$; $k(i) =1$ \\ for $i \in \{1, \ldots , 6\}$ where 
$i \ne j$ and $i \ne \ell $} \end{array} \right. 
\end{displaymath} 
is the set of roots $\mathfrak{R}$ whose root vectors $X_{\widehat{0}}; X_{\widehat{\alpha }}, 
\widehat{\alpha } \in \widehat{\mathcal{R}}$ in $\underline{\mathfrak{g}}$ span a complex vector space 
$\widetilde{\mathfrak{n}}$. \medskip 

We now prove some properties of $\widehat{\mathcal{R}}$, which we use in claim A2 to show that 
$\widetilde{\mathfrak{n}}$ is a sandwich algeba whose center is spanned by $X_{\widehat{0}}$. 
First, for every $\widehat{\alpha } \in \widehat{\mathcal{R}}$
\begin{equation}
\widehat{0} + \widehat{\alpha } \notin \widehat{\mathcal{R}}, 
\label{eq-Aonec}
\end{equation}
because $-{\varepsilon}_7 + \widehat{\alpha} \notin \underline{\mathcal{R}}$. Second, by definition of 
${\widehat{\alpha}}_{jk}$ and ${\widehat{\alpha}}^{jk}$ we have 
\begin{equation}
{\widehat{\alpha}}_{jk} +{\widehat{\alpha}}^{jk} = \widehat{0}
\label{eq-Atwoc}
\end{equation}
for every $j$, $k \in \{ 1, 2, \ldots , 6 \}$ with $j < k$. Equation (\ref{eq-Atwoc}) also holds when 
$(j,k) = (1,0)$. Third 
\begin{equation}
{\widehat{\alpha }}_{jk} + {\widehat{\alpha }}^{\ell m} \notin \widehat{\mathcal{R}}, 
\, \, \, \mbox{when $(\ell ,m) \ne (j,k)$}
\label{eq-Athreec}
\end{equation}
because at least one of the first six components of ${\widehat{\alpha }}_{jk} + {\widehat{\alpha }}^{\ell m}$ 
is equal to $0$ or $\pm 1$ and its seventh component is equal to $-1$. Thus the set of roots $\widehat{\mathcal{R}}$ is closed under addition. Using claim A1 it follows that $\widetilde{\mathfrak{n}}$ is the nilradical of the Lie subalgebra $\widetilde{\mathfrak{g}} = \mathfrak{g} \oplus \widetilde{\mathfrak{n}}$. Moreover, $\widetilde{\mathfrak{n}}$ is an ideal of $\underline{\mathfrak{g}}$. \bigskip

We now prove \medskip 

\noindent \textbf{Claim A1} Suppose that the set of roots $\widehat{\mathcal{R}}$ is closed under addition. 
Then the space $\widetilde{\mathfrak{n}}$ spanned by the root vectors $X_{\widehat{0}}$ and 
$X_{\widehat{\alpha }}$, $\widehat{\alpha } \in \widehat{\mathcal{R}}$ in 
the simple Lie algebra $\underline{\mathfrak{g}}$ is the nilradical of the Lie subalgebra 
$\widetilde{\mathfrak{g}} = \mathfrak{g} \oplus \widetilde{\mathfrak{n}}$ and is an ideal of 
$\underline{\mathfrak{g}}$. \medskip 

\noindent \textbf{Proof.} By construction ${\spann }_{\C} \{ X_{\widehat{0}} ;  
X_{\widehat{\alpha }}, \, \widehat{\alpha } \in \widehat{\mathcal{R}} \} $ is a complex vector subspace 
of $\underline{\mathfrak{g}}$, which is invariant under ${\ad }_{\mathfrak{h}}$, where 
$\mathfrak{h}$ is a Cartan subalgebra of $\mathfrak{g}$. By definition ${\mathcal{R}}^{-} \subseteq 
\underline{\mathcal{R}}$ is closed under addition. Hence $\widetilde{\mathfrak{n}}$ is 
a nilpotent subalgebra of $\underline{\mathfrak{g}}$. Since $\mathfrak{R} \cup \widehat{\mathfrak{R}} 
\subseteq \underline{\mathcal{R}}$ is closed under addition, $\widetilde{\mathfrak{g}} = 
\mathfrak{g} \oplus \widetilde{\mathfrak{n}}$ is a Lie subalgebra of $\underline{\mathfrak{g}}$ and 
$\widetilde{\mathfrak{n}}$ is an ideal, which is the nilradical of $\widetilde{\mathfrak{g}}$, because 
$\mathfrak{g}$ is a simple Lie algebra. \hfill $\square $ \medskip 

\noindent \textbf{Claim 2A} Suppose that $\widehat{\alpha } \in \widehat{\mathcal{R}}$ satisfies 
the following three conditions 
\begin{equation}
\widehat{\alpha}  + \widehat{0} \notin \widehat{\mathcal{R}}; 
\label{eq-Aonestar} 
\end{equation} 
for every $\widehat{\alpha } \in \widehat{\mathcal{R}}$ there is a 
$\widehat{\beta } \in \widehat{\mathcal{R}}$ such that 
\begin{equation} 
\widehat{\alpha } + \widehat{\beta } = \widehat{0};
\label{eq-Atwostar} 
\end{equation}
for every $\widehat{\alpha }$ and $\widehat{\beta } \in 
\widehat{\mathcal{R}}$ with $\widehat{\alpha } + \widehat{\beta } \ne \widehat{0}$ one has 
\begin{equation}
\widehat{\alpha } + \widehat{\beta } \notin \widehat{\mathcal{R}} .
\label{eq-Athreestar}
\end{equation} 
Then $\widetilde{\mathfrak{n}}$ is a sandwich algebra whose center is spanned 
by $X_{\widehat{0}}$. Moreover $\widetilde{\mathfrak{n}}$ is a Heisenberg algebra. \medskip 

\noindent \textbf{Proof.} Using equation (\ref{eq-Aonestar}) we obtain $[X_{\widehat{\alpha}} , X_{\widehat{0}} ] 
=0 $ for all $\widehat{\alpha } \in \widehat{\mathcal{R}}$. But $[X_{\widehat{0}}, X_{\widehat{0}}] = 0$. 
So $X_{\widehat{0}}$ spans the center of $\widetilde{\mathfrak{n}}$. From equation 
(\ref{eq-Atwostar}) it follows that $[X_{\widehat{\alpha }}, X_{\widehat{\beta }}] = X_{\widehat{0}}$. So the Lie subalgebra $\widetilde{\mathfrak{n}}$ is not abelian. Using equations (\ref{eq-Aonestar}), 
(\ref{eq-Atwostar}) and (\ref{eq-Athreestar}) we see that 
no sum of three roots in $\widehat{\mathcal{R}}$ lies in $\mathcal{R}$. Thus the 
bracket $[ \cdot , [ \cdot , \cdot ]]$ of their corresponding root vectors vanishes. So 
$\widetilde{\mathfrak{n}}$ is a sandwich. From equations (\ref{eq-Aonestar}), (\ref{eq-Atwostar}) and 
(\ref{eq-Athreestar}) we see that the vectors in $\widetilde{\mathfrak{n}}$ satisfy the bracket relations for a Heisenberg algebra. \hfill $\square $ 

\label{references}

\end{document}